%Paper: alg-geom/9210008
%From: David Cox <dac@cs.amherst.edu>
%Date: Thu, 22 Oct 92 13:50:23 EDT
%Date (revised): Mon, 21 Jun 93 14:16:47 EDT

% Here is the TeX file for the revised paper "The Homogeneous Coordinate Ring
% of a Toric Variety'' by David A. Cox.

% The file is in pure TeX (C Version 3.0), but it requires the use of
% the AMS fonts (which are invoked in a non-standard way).

% Thus there are two parts to the file:

% 1) MACROS: The macros used in the main part of the paper.
% 2) TEXT: The text of the paper.

%%%%%%%%%%%%%%%%%%%%%%%%%%%%%%%%%%%%%%%%%%%%%%%%%%%%%%%%%%%%%%%%%%%%%%

% 1) MACROS: Here are the macros used in "The Homogeneous Coordinate
%            Ring of a Toric Variety''.

% The page format:

\magnification=\magstep1
\baselineskip=14pt
\abovedisplayskip=10pt plus 2pt minus 2pt
\belowdisplayskip=10pt plus 2pt minus 2pt

% Special macro for math display:

\def\d{$$}

% Special macros:

\font\bigbf=cmbx10 scaled \magstep1 % large font
\def\ref{\hangindent=\parindent \hangafter=1 \noindent} % used for references

% Some commonly used greek letters:

\def\a{\alpha}

% Blackboard Roman letters

\font\tenmsy=msbm10
\font\sevenmsy=msbm7
\font\fivemsy=msbm5
%\font\tenmsy=msym10
%\font\sevenmsy=msym7
%\font\fivemsy=msym5
\def\bb#1{\mathchoice{\hbox{\tenmsy #1}}{\hbox{\tenmsy #1}}
      	{\hbox{\sevenmsy #1}}{\hbox{\fivemsy #1}}}

\def\P{{\bb P}}

\def\z{{\bb Z}}
\def\q{{\bb Q}}
\def\r{{\bb R}}
\def\cx{{\bb C}}  % \c is already in use

% Other formatting conventions

\def\prf{\noindent {\bf Proof.\ \ }}
\def\square{\smash{\raise8pt\hbox{\leavevmode\hbox{\vrule\vtop{\vbox
	{\hrule width9pt}\kern9pt\hrule width9pt}\vrule}}}}
\def\qed{\hfill\square}

% Other math symbols

\font\tenmsx=msam10

%\font\tenmsx=msxm10
%\font\sevenmsx=msxm7
%\font\fivemsx=msxm5

\newfam\msx

\textfont8=\tenmsy
\textfont9=\tenmsx

\mathchardef\subsetne = "3824
\mathchardef\supsetne = "3825
\mathchardef\notdiv = "282D
\mathchardef\semidir = "286F
\mathchardef\ontoarrow = "3910

\def\mapname#1{\smash{\mathop{\longrightarrow}\limits^{#1}}}

%%%%%%%%%%%%%%%%%%%%%%%%%%%%%%%%%%%%%%%%%%%%%%%%%%%%%%%%%%%%%%%%%%%%%%%%

% 2) TEXT: Here is the text of the revised paper "The Homogeneous Coordinate
%          Ring of a Toric Variety''.

\centerline{\bigbf The Homogeneous Coordinate Ring of a Toric Variety}
\centerline{\it Revised June 21, 1993}
\smallskip
\centerline{David A. Cox}
\centerline{Department of Mathematics and Computer Science}
\centerline{Amherst College}
\centerline{Amherst, MA 01002}
\medskip

	This paper will introduce the {\it homogeneous coordinate
ring\/} $S$ of a toric variety $X$.  The ring $S$ is a polynomial ring
with one variable for each one-dimensional cone in the fan $\Delta$
determining $X$, and $S$ has a natural grading determined by the
monoid of effective divisor classes in the Chow group $A_{n-1}(X)$ of
$X$ (where $n = \dim X$).  Using this graded ring, we will show that
$X$ behaves like projective space in many ways.  The paper is
organized into four sections as follows.

	In \S1, we define the homogeneous coordinate ring $S$ of $X$
and compute its graded pieces in terms of global sections of certain
coherent sheaves on $X$.  We also define a monomial ideal $B \subset
S$ that describes the combinatorial structure of the fan $\Delta$.  In
the case of projective space, the ring $S$ is just the usual
homogeneous coordinate ring $\cx[x_0,\dots,x_n]$, and the ideal $B$ is
the ``irrelevant'' ideal $\langle x_0,\dots,x_n\rangle$.

	Projective space $\P^n$ can be constructed as the quotient
$(\cx^{n+1} - \{0\})/\cx^*$.  In \S2, we will see that there is a
similar construction for any toric variety $X$.  In this case, the
algebraic group $G = {\rm Hom}_\z(A_{n-1}(X),\cx^*)$ acts on an affine
space $\cx^{\Delta(1)}$ such that the categorical quotient
$(\cx^{\Delta(1)} - Z)/G$ exists and is isomorphic to $X$.  The
exceptional set $Z$ is the zero set of the ideal $B$ defined in
\S1.  If $X$ is simplicial (meaning that the fan $\Delta$ is
simplicial), then $X \simeq (\cx^{\Delta(1)}-Z)/G$ is a geometric
quotient, so that elements of $\cx^{\Delta(1)} - Z$ can be
regarded as ``homogeneous coordinates'' for points of $X$.

	On $\P^n$, there is a correspondence between sheaves and
graded $\cx[x_0,\dots,x_n]$-modules.  For any toric variety, we will
see in \S3 that finitely generated graded $S$ modules give rise to a
coherent sheaves on $X$, and when $X$ is simplicial, every coherent
sheaf arises in this way.  In particular, every closed subscheme of
$X$ is determined by a graded ideal of $S$.  We will also study
the extent to which this correspondence fails to be is one-to-one.

	Another feature of $\P^n$ is that the action of $PGL(n+1,\cx)$
on $\P^n$ lifts to an action of $GL(n+1,\cx)$ on $\cx^{n+1}-\{0\}$.
In \S4, we will see that for any complete simplicial toric variety
$X$, the action of ${\rm Aut}(X)$ on $X$ lifts to an action of an
algebraic group $\widetilde{{\rm Aut}}(X)$ on
$\cx^{\Delta(1)}-Z$.  This group is a extension of ${\rm Aut}(X)$
by the group $G$ defined above.  We will give an explicit description
of $\widetilde{\rm Aut}(X)$, and in particular, we will see that the
roots of ${\rm Aut}(X)$ (as defined in [Demazure]) have an especially
nice interpretation.

	For simplicity, we will work over the complex numbers $\cx$.
Our notation will be similar to that used by [Fulton] and [Oda].  I
would like to thank Bernd Sturmfels for stimulating my interest in
toric varieties.

\bigskip
\noindent{\bf \S1. The Homogeneous Coordinate Ring}
\bigskip

	Let $X$ be the toric variety determined by a fan $\Delta$ in
$N \simeq \z^n$.  As usual, $M$ will denote the $\z$-dual of $N$, and
cones in $\Delta$ will be denoted by $\sigma$.  The one-dimensional
cones of $\Delta$ form the set $\Delta(1)$, and given $\rho \in
\Delta(1)$, we let $n_\rho$ denote the unique generator of $\rho\cap
N$.  If $\sigma$ is any cone in $\Delta$, then $\sigma(1) = \{\rho \in
\Delta(1): \rho \subset \sigma\}$ is the set of one-dimensional faces
of $\sigma$.  We will assume that $\Delta(1)$ spans $N_\r =
N\otimes_\z \r$.

	Each $\rho \in \Delta(1)$ corresponds to an irreducible
$T$-invariant Weil divisor $D_\rho$ in $X$, where $T = N\otimes_\z
\cx^*$ is the torus acting on $X$.  The free abelian group of
$T$-invariant Weil divisors on $X$ will be denoted $\z^{\Delta(1)}$.
Thus an element $D \in \z^{\Delta(1)}$ is a sum $D = \sum_{\rho}
a_\rho D_\rho$.  The $T$-invariant Cartier divisors form a subgroup
${\rm Div}_{T}(X) \subset \z^{\Delta(1)}$.

	Each $m \in M$ gives a character $\chi^m : T \to \cx^*$, and
hence $\chi^m$ is a rational function on $X$.  As is well-known,
$\chi^m$ gives the Cartier divisor ${\rm div}(\chi^m) = - \sum_{\rho}
\langle m,n_\rho\rangle D_\rho$.  We find it convenient to ignore
the minus sign, so that we will consider the map
\d
M \longrightarrow \z^{\Delta(1)}\ \hbox{defined by}\ m \mapsto
D_m = \textstyle{\sum_\rho \langle m,n_\rho\rangle D_\rho}\ .
\d
This map is injective since $\Delta(1)$ spans $N_\r$.  By [Fulton,
\S3.4], we have a commutative diagram
\d
\matrix{0 & \longrightarrow & M & \longrightarrow & {\rm Div_{T}}(X)
	& \longrightarrow & {\rm Pic}(X) & \longrightarrow & 0\cr
& & \| & & \downarrow & & \downarrow & & \cr
0 & \longrightarrow & M & \longrightarrow & \z^{\Delta(1)} &
	\longrightarrow & A_{n-1}(X) & \longrightarrow &0\cr}
	\leqno(1)
\d
where the rows are exact and the vertical arrows are inclusions.  Thus
a divisor $D \in \z^{\Delta(1)}$ determines an element $\a = [D] \in
A_{n-1}(X)$.  Note that $n$ is the dimension of $X$.

	For each $\rho \in \Delta(1)$, introduce a variable $x_\rho$,
and consider the polynomial ring
\d
S = \cx[x_\rho : \rho \in \Delta(1)]\ .
\d
We will usually write this as $S = \cx[x_\rho]$.  Note that a monomial
$\prod_\rho x_\rho^{a_\rho}$ determines a divisor $D = \sum_\rho
a_\rho D_\rho$, and to emphasize this relationship, we will write the
monomial as $x^D$.  We will grade $S$ as follows:
\d
\hbox{the {\it degree\/} of a monomial $x^D \in S$ is ${\rm deg}(x^D)
	= [D] \in A_{n-1}(X)$\ . }
\d
Using the exact sequence (1), it follows that two monomials
$\prod_\rho x_\rho^{a_\rho}$ and $\prod_\rho x_\rho^{b_\rho}$ in $S$
have the degree if and only if there is some $m \in M$ such that
$a_\rho = \langle m,n_\rho\rangle + b_\rho$ for all $\rho$.   Then let
\d
S_\a = \bigoplus_{{\rm deg}(x^D) = \a} \cx\cdot x^D\ ,
\d
so that the ring $S$ can be written as the direct sum
\d
S = \bigoplus_{\a \in A_{n-1}(X)} S_\a\ .
\d
Note also that $S_\a\cdot S_\beta \subset S_{\a+\beta}$.  We call $S$
the {\it homogeneous coordinate ring\/} of the toric variety $X$.  Of
course, ``homogeneous'' means with respect to the above grading.  We
should mention that the ring $S$, without the grading, appears in the
Danilov-Jurkiewicz description of the cohomology of a simplicial toric
variety (see [Danilov, \S10], [Fulton, \S5.2] or [Oda, \S3.3]).

	Here are some examples of what the ring $S$ looks like:
\medskip

\noindent {\bf Projective Space:} When $X = \P^n$, it is easy to check
that $S$ is usual homogeneous coordinate ring $\cx[x_0,\dots,x_n]$
with the standard grading.
\medskip

\noindent{\bf Weighted Projective Space:} When $X = \P(q_0,\dots,q_n)$,
then $S$ is the ring $\cx[x_0,\dots,x_n]$, where the grading is
determined by giving the variable $x_i$ weight $q_i$.
\medskip

\noindent{\bf Product of Projective Spaces:} When $X = \P^n\times\P^m$,
then $S = \cx[x_0,\dots,x_n;y_0,\dots,y_m]$.  Here, the grading is the
usual bigrading, where a polynomial has bidegree $(a,b)$ if it is
homogeneous of degree $a$ (resp. $b$) in the $x_i$ (resp. $y_j$).
\medskip

	Our first result about the ring $S$ shows that the graded
pieces of $S$ are isomorphic to the global sections of certain sheaves
on $X$:

\proclaim Proposition 1.1.
\vskip0pt
\item{(i)} If $\a = [D] \in A_{n-1}(X)$, then there is an
isomorphism
\d
\phi_D: S_\a \simeq H^0(X,{\cal O}_X(D))\ ,
\d
where ${\cal O}_X(D)$ is the coherent sheaf on X determined by the
Weil divisor $D$ (see [Fulton, \S3.4]).
\item{(ii)} If $\a =[D]$ and $\beta = [E]$, then there is a
commutative diagram
\d
\matrix{S_\a\otimes S_\beta & \longrightarrow & S_{\a+\beta}\cr
	\downarrow && \downarrow \cr
	H^0(X,{\cal O}_X(D))\otimes H^0(X,{\cal O}_X(E)) &
		\longrightarrow & H^0(X,{\cal O}_X(D+E))\cr}
\d
where the top arrows is multiplication, the bottom arrow is tensor
product, and the vertical arrows are the isomorphisms
$\phi_D\otimes\phi_E$ and $\phi_{D+E}$.

\prf Suppose that $D = \sum_\rho a_\rho D_\rho$.  By [Fulton,
\S3.4], we know that $H^0(X,{\cal O}_X(D)) = \oplus_{\a \in P_D\cap M}
\cx\cdot\chi^m$, where
\d
	P_D = \{m \in M_\r : \langle m,n_\rho\rangle \ge -a_\rho\
		\hbox{for all}\ \rho\}\ .
\d
Given $m \in P_D\cap M$, let $D_m = \sum_\rho \langle m,n_\rho\rangle
D_\rho$, so that $x^{D_m + D} = \prod_\rho x_\rho^{\langle
m,n_\rho\rangle + a_\rho}$.  Then $m \mapsto x^{D_m+D}$ defines a map
from $P_D\cap M$ to the monomials in $S_\a$ (the monomial is in $S$
since $m \in P_D$, and it has degree $\a$ since $m \in M$).  This map
is one-to-one since $M \to \z^{\Delta(1)}$ is injective, and it is
easy to see that every monomial in $S_\a$ arises in this way.  Then
$\phi_D(x^{D_m+D}) = \chi^m$ gives the desired isomorphism.  The proof
of the second part of the proposition is straightforward and is
omitted.  \qed
\medskip

	We get the following corollary (see [Fulton, \S3.4]):

\proclaim Corollary 1.2. Let $X$ be a complete toric variety.  Then:
\vskip0pt
\item{(i)} $S_\a$ is finite dimensional for every $\a$, and in
particular, $S_0 = \cx$.
\item{(ii)} If $\a = [D]$ for an effective divisor $D = \sum_\rho
a_\rho D_\rho$, then $\dim_\cx S_\a = |P_D\cap M|$.
\vskip0pt

	In the complete case, the graded ring $S$ has some additional
structure.  Namely, there is a natural order relation on the monomials
of $S$ defined as follows: if $x^D, x^E \in S$, then set
\d
x^D < x^E \iff \hbox{there is $x^F \in S$ with $x^D |
x^F$, $x^D \ne x^F$ and $\deg(x^F) = \deg(x^E)$}\ .
\d
This relation has the following properties:

\proclaim Lemma 1.3. When $X$ is a complete toric variety, the order
relation defined above is transitive, antisymmetric and multiplicative
(meaning that $x^D < x^E \Rightarrow x^{D+F} < x^{E+F}$) on the
monomials of $S$.

\prf The transitive and multiplicative properties are trivial to
verify.  To prove antisymmety, note that $x^D | x^F$ and $\deg(x^F) =
\deg(x^E)$ implies that $[E]-[D] = [F] - [D] = [F-D]$ is the class of
an effective divisor.  Thus, if $x^D < x^E$ and $x^E < x^D$, then $[E]
- [D]$ and its negative are effective classes.  Since $X$ is complete,
irreducible and reduced, we must have $[D] = [E]$.  Then, for $x^F$ as
above, we would have $x^{F-D} \in S_0 \simeq \cx$, which would imply
$x^F = x^D$.  This contradicts the definition of $x^D < x^E$, and
antisymmetry is proved.\qed
\medskip

	In the case of $\P^n$, this gives the usual ordering by total
degree.  For a complete simplicial toric variety $X$, we will use the
order relation of Lemma 1.3 in \S4 when we study the automorphism
group of $X$.

	An important observation is that the theory developed so far
depends only on the one-dimensional cones $\Delta(1)$ of the fan
$\Delta$.  More precisely, if $\Delta$ and $\Delta'$ are fans in $N$
with $\Delta(1) = \Delta'(1)$, then $A_{n-1}(X) = A_{n-1}(X')$, and
the corresponding rings $S$ and $S'$ are also equal as graded rings.
It follows that we cannot reconstruct the fan from $S$ and its
grading---more information is needed.  The ring tells us about
$\Delta(1)$, but we need something else to tell us which cones belong
in $\Delta$.

	The crucial object is the following ideal of $S$.  For a cone
$\sigma \in \Delta$, let $\hat\sigma$ be the divisor $\sum_{\rho
\notin \sigma(1)} D_\rho$, and let the corresponding monomial in $S$ be
$x^{\hat\sigma} = \prod_{\rho \notin \sigma(1)} x_\rho$.  Then define
$B = B_\Delta \subset S$ to be the ideal generated by the
$x^{\hat\sigma}$, i.e.,
\d
B = \langle x^{\hat\sigma} : \sigma \in \Delta \rangle \subset S\ .
\d
Note the $B$ is in fact generated by the $x^{\hat\sigma}$ as $\sigma$
ranges over the {\it maximal\/} cones of $\Delta$.  Also, if $\Delta$
and $\Delta'$ are fans in $N$ with $\Delta(1) = \Delta'(1)$, then
$\Delta = \Delta'$ if and only if the corresponding ideals $B, B'
\subset S$ satisfy $B = B'$.  This follows because $\{x^{\hat\sigma} :
\hbox{$\sigma$ is a maximal cone of $\Delta$}\}$ is the unique minimal
basis of the monomial ideal $B$.

	When $X$ is a projective space or weighted projective space,
the ideal $B$ is just the ``irrelevant'' ideal $\langle x_0,\dots,x_n
\rangle \subset \cx[x_0,\dots,x_n]$.  A basic theme of this paper is
that the pair $B \subset S$ plays the role for an arbitrary toric
variety that $\langle x_0,\dots,x_n\rangle \subset \cx[x_0,\dots,x_n]$
plays for projective space.  In particular, we can think of the
variety
\d
Z = {\bf V}(B) = \{ x \in \cx^{\Delta(1)} : x^{\hat\sigma} = 0\
\hbox{for all $\sigma \in \Delta$}\} \subset \cx^{\Delta(1)}
\d
as the ``exceptional'' subset of $\cx^{\Delta(1)}$.  In \S2, we will
show how to construct the toric variety $X$ from $\cx^{\Delta(1)} -
Z$.

\proclaim Lemma 1.4. $Z \subset \cx^{\Delta(1)}$ has codimension
at least two.

\prf Since $B' = \langle x^{\hat\rho} : \rho \in \Delta(1)\rangle
\subset B$, we have $Z = {\bf V}(B) \subset {\bf V}(B')$.  But
${\bf V}(B')$ is the union of all codimension two coordinate
subspaces, and the lemma follows. \qed
\bigskip

\noindent{\bf \S2. Homogeneous Coordinates}
\bigskip

	Given a toric variety $X$, the Chow group $A_{n-1}(X)$ defined
in \S1 is a finitely generated abelian group of rank $d - n$, where $d
= |\Delta(1)|$.  It follows that $G = {\rm Hom}_\z(A_{n-1}(X),\cx^*)$
is isomorphic to a product of the torus $(\cx^*)^{d-n}$ and the finite
group ${\rm Hom}_\z(A_{n-1}(X)_{tor}, \q/\z)$.  Furthermore, if we
apply ${\rm Hom}(-,\cx^*)$ to the bottom exact sequence of (1), then
we get the exact sequence
\d
1 \longrightarrow G \longrightarrow (\cx^*)^{\Delta(1)}
	\longrightarrow T \longrightarrow 1\ ,
\d
where $T = N\otimes\cx^* = {\rm Hom}_\z(M,\cx^*)$ is the torus acting
on $X$.  Note also that $A_{n-1}(X)$ is naturally isomorphic to the
character group of $G$.

	Since $(\cx^*)^{\Delta(1)}$ acts naturally on
$\cx^{\Delta(1)}$, the subgroup $G \subset (\cx^*)^{\Delta(1)}$ acts
on $\cx^{\Delta(1)}$ by
\d
g\cdot t = (g([D_\rho])\,t_\rho)
\d
for $g : A_{n-1}(X) \to \cx^*$ in $G$ and $t = (t_\rho)$ in
$\cx^{\Delta(1)}$.  We get a corresponding representation of $G$ on $S
= \cx[x_\rho]$, and given a character $\a \in A_{n-1}(X)$ of $G$, the
graded piece $S_\a$ is the $\a$-eigenspace of this representation.

	Using the action of $G$ on $\cx^{\Delta(1)}$, we can describe
the toric variety $X$ as follows:

\proclaim Theorem 2.1.  Let $X$ be the toric variety determined by the
fan $\Delta$, and let $Z = {\bf V}(B) \subset \cx^{\Delta(1)}$ be
as in \S1.  Then:
\vskip0pt
\item{(i)} The set $\cx^{\Delta(1)} - Z$ is invariant under the
action of the group $G$.
\item{(ii)} $X$ is naturally isomorphic to the categorical quotient of
$\cx^{\Delta(1)} - Z$ by $G$.
\item{(iii)} $X$ is the geometric quotient of $\cx^{\Delta(1)} - Z$ by
$G$ if and only if $X$ is simplicial (i.e., the fan $\Delta$ is
simplicial).
\vskip0pt

	In the analytic category, part (iii) of Theorem 2.1 was known
to M.~Audin (see Chapter VI of [Audin], which mentions earlier work by
Delzant and Kirwan), and versions of the theorem were known to
V.~Batyrev (unpublished) and J.~Fine (just now appearing in [Fine]).
As this paper was being written, I.~Musson discovered a result closely
related to part (ii) of the theorem (see [Musson]), and more recently
he proved an analog of part (iii).

	In the first version of this paper, part (iii) only stated
that the geometric quotient exists when $X$ is simplicial.  The
referee asked if the converse were true, and simultaneously Musson
proved this in a slightly different content.  His argument is used
below.
\medskip

\noindent {\bf Proof of Theorem 2.1.}  Given a cone $\sigma \in
\Delta$, let
\d
U_\sigma = \cx^{\Delta(1)} - {\bf V}(x^{\hat\sigma}) = \{ x \in
\cx^{\Delta(1)} : x^{\hat\sigma} \ne 0\}\ .
\d
Note that $U_\sigma$ is a $G$-invariant affine open subset of
$\cx^{\Delta(1)}$ and that
\d
\cx^{\Delta(1)} - Z = \bigcup_{\sigma \in \Delta} U_\sigma\ .
\d
This proves part (i) of the theorem.

	For parts (ii) and (iii) of the theorem, we first need to
study the quotient of $G$ acting on the affine variety $U_\sigma$.
The coordinate ring of $U_\sigma$ is the localization of $S$ at
$x^{\hat\sigma}$, which we will denote $S_\sigma $.  Note that
$S_\sigma$ has a natural grading by $A_{n-1}(X)$.

	The following lemma describes the degree zero part of
$S_\sigma$:

\proclaim Lemma 2.2. Let $\sigma \in \Delta$ be a cone and let
$\check\sigma \subset M_\r$ be its dual cone.  Then there is a natural
isomorphism of rings $\cx[\check\sigma\cap M] \simeq (S_\sigma)_0$.

\prf Given $m \in M$, let $x^{D_m} = \prod_\rho x_\rho^{\langle
m,n_\rho\rangle}$, and note that
\d
x^{D_m} \in S_\sigma \iff \langle m,n_\rho\rangle \ge 0\
	\hbox{for all $\rho \in \sigma(1)$}
	\iff m \in \check\sigma\ .
\d
Since a monomial has degree zero if and only if it is of the form
$x^{D_m}$ for some $m \in M$, we get a one-to-one
correspondence between $\check\sigma\cap M$ and monomials in
$(S_\sigma)_0$.  Then the map sending $m \in \check\sigma\cap M$ to
$x^{D_m} \in (S_\sigma)_0$ gives the desired ring isomorphism.\qed
\medskip

	The action of $G$ on $U_\sigma$ induces an action on
$S_\sigma$, and as we observed earlier about $S$, the graded pieces
are the eigenspaces for the characters of $G$.  In particular, the
invariants of $G$ acting on $S_\sigma$ are precisely the polynomials
of degree zero.  Thus Lemma 2.2 implies
\d
(S_\sigma)^G = (S_\sigma)_0 \simeq \cx[\check\sigma\cap M]\ .
\d
Since $G$ is a reductive group acting on the affine scheme $U_\sigma$,
standard results from invariant theory (see [Fogarty and Mumford,
Theorem 1.1]) imply that the categorical quotient is given by
\d
U_\sigma/G = {\rm Spec}((S_\sigma)^G) = {\rm Spec}(\cx[\check\sigma\cap
	M])\ .
\d
Thus $U_\sigma/G$ is the affine toric variety $X_\sigma$ determined by
the cone $\sigma$.  Note that $X_\sigma$ is an affine piece of the
toric variety $X$.

	The next step is to see how the quotients $U_\sigma/G$ fit
together as we vary $\sigma$.  So let $\tau$ be a face of a cone
$\sigma \in \Delta$.  It is well-known that $\tau =
\sigma\cap\{m\}^\perp$ for some $m \in \check\sigma\cap M$.  Then we
can relate $(S_\tau)_0$ and $(S_\sigma)_0$ as follows:

\proclaim Lemma 2.3. Let $\tau$ be a face of $\sigma$, and let $m \in
\check\sigma\cap M$ be as above.  Under the isomorphism of Lemma
2.2, $m$ maps to the monomial $x^{D_m} = \prod_\rho x_\rho^{\langle
m,n_\rho\rangle}\in (S_\sigma)_0$.  Then:
\vskip0pt
\item{(i)} There is a natural identification
\d
(S_\sigma)_0 = \big((S_\tau)_0\big)_{x^{D_m}}\ .
\d
\item{(ii)} There is a commutative diagram
\d
\matrix{(S_\tau)_0 & \longrightarrow & \ \ (S_\sigma)_0 =
	\big((S_\tau)_0\big)_{x^{D_m}} \cr
\downarrow && \quad\quad\downarrow\hfill \cr
\cx[\check\tau\cap M] & \longrightarrow & \cx[\check\sigma\cap M] =
	\cx[\check\tau\cap M]_m\cr}
\d
where the horizontal maps are localization maps and the vertical
arrows are the isomorphisms of Lemma 2.2.
\vskip0pt

\prf Since $\tau \in \sigma\cap\{m\}^\perp$, $x^{D_m}$ has positive
exponent for $x_\rho$ when $\rho \in \sigma(1) - \tau(1)$ and zero
exponent when $\rho \in \tau(1)$.  It follows immediately that $S_\tau
= (S_\sigma)_{x^{D_m}}$.  Part (i) of the lemma follows since
localization commutes with taking elements of degree zero (because
$x^{D_m}$ already has degree zero).  It is trivial to check that the
diagram of part (ii) commutes, and the lemma is proved.\qed
\medskip

	We can now prove part (ii) of Theorem 2.1.  Namely, Lemma 2.2
gives us categorical quotients $U_\sigma/G \cong X_\sigma$, and by the
compatibility of Lemma 2.3, these patch to give a categorical quotient
$(\cx^{\Delta(1)} - Z)/G \cong X$.

	For part (iii) of the theorem, first assume that $X$ is
simplicial.  To show that $X$ is a geometric quotient, it suffices to
show that each $X_\sigma$ is the geometic quotient of $G$ acting on
$U_\sigma$.  By [Fogarty and Mumford, Amplification 1.3], we need only
show that the orbits of $G$ acting on $U_\sigma$ are closed.

	Fix a point $t = (t_\rho)$ of $U_\sigma$.  In Lemma 2.2 we saw
that $m \in \check\sigma\cap M$ gives the monomial $x^{D_m} =
\prod_\rho x_\rho^{\langle m,n_\rho\rangle} \in (S_\sigma)_0$.  Then
consider the following subvariety $V \subset U_\sigma$:
\d
V = \{x = (x_\rho) \in U_\sigma : x^{D_m} = t^{D_m}\ \hbox{for all
$m \in \check\sigma\cap M$}\}\ .
\d
We will prove that $V = G\cdot t$.  First observe that $t \in V$, and
since $x^{D_m}$ is $G$-invariant (it has degree zero), it follows that
$G\cdot t \subset V$.

	To prove the other inclusion, let $u = (u_\rho) \in V$.  We
must find $g:A_{n-1}(X) \to \cx^*$ such that $u_\rho =
g([D_\rho])\,t_\rho$ for all $\rho \in \Delta(1)$.  The first step is
to show that
\d
\hbox{for all $\rho \in \Delta(1)$},\ u_\rho \ne 0 \iff t_\rho \ne 0
	\ .\leqno{(2)}
\d
Note that (2) is trivially true when $\rho \notin \sigma(1)$.  So
suppose that $\rho \in \sigma(1)$.  Since $X$ is simplicial, the set
$\{n_{\rho'}: \rho' \in \sigma(1)\}$ is linearly independent.  Hence
we can find $m \in M$ such that
\d
\langle m,n_\rho\rangle > 0\ {\rm and}\ \langle m,n_{\rho'}\rangle = 0\
	\hbox{for all $\rho' \in \sigma(1) - \{\rho\}$}\ .
\d
Setting $a = \langle m,n_\rho\rangle > 0$, we have
\d
\eqalign{u^{D_m} &= \prod_{\rho' \notin \sigma(1)} u_{\rho'}^{\langle
	m,n_{\rho'}\rangle} \cdot u_\rho^a\cr
t^{D_m} &= \prod_{\rho' \notin \sigma(1)} t_{\rho'}^{\langle
	m,n_{\rho'}\rangle} \cdot t_\rho^a\ .\cr}
\d
Since these are equal, and since $u_{\rho'}$ and $t_{\rho'}$ are
nonzero for $\rho' \notin \sigma(1)$, it follows that $u_\rho \ne 0$
if and only if $t_\rho \ne 0$.  This completes the proof of (2).

	Now let $A = \{\rho : t_\rho = 0\} = \{\rho : u_\rho = 0\}$
and note that $A \subset \sigma(1)$.  For $\rho \in \Delta(1) - A$,
set
\d
g_A(D_\rho) = {u_\rho \over t_\rho}\ .
\d
This defines a map $g_A : \z^{\Delta(1) - A} \to \cx^*$.  To relate
this to $G = {\rm Hom}_\z(A_{n-1}(X),\cx^*)$, let $i_A : M \to \z^{A}$
be the composition $M \to \z^{\Delta(1)} \to \z^{A}$, where the last
map is projection.  Then we have the commutative diagram
\d
\matrix{&&&&0&&&&\cr
&&&&\downarrow&&&&\cr
0 & \longrightarrow & 0 & \longrightarrow & M & = & M\hfill &
	\longrightarrow & 0\cr
&& \downarrow && \downarrow && \,\downarrow\ i_A &&\cr
0 & \longrightarrow & \z^{\Delta(1)-A} & \longrightarrow &
	\z^{\Delta(1)} & \longrightarrow & \z^{A}\hfill &
	\longrightarrow & 0\cr
&& \| && \downarrow &&&&\cr
&&\z^{\Delta(1)-A} & \longrightarrow & A_{n-1}(X) &&&&\cr
&&\downarrow&&\downarrow&&&&\cr
&&0&&0&&&&\cr}
\d
Since the columns and two middle rows are exact, the snake lemma gives
us an exact sequence
\d
{\ker}(i_A) \longrightarrow \z^{\Delta(1)-A} \longrightarrow
	A_{n-1}(X)\ .\leqno(3)
\d

	We next claim that $g_A(\ker(i_A)) = \{1\}$.  An easy diagram
chase shows that for $m \in \ker(i_A)$, we have
\d
g_A(m) = {u^{D_m} \over t^{D_m}}\ .
\d
When $m \in \check\sigma\cap M$, we know that $u^{D_m} = t^{D_m}$, so
that $g_A(m) = 1$ whenever $m \in \check\sigma\cap\ker(i_A)$.  If we
can show that this semigroup generates the group $\ker(i_A)$, then our
claim will follow.  But $\{n_\rho : \rho \in \sigma(1)\}$ is linearly
independent since $\sigma$ is simplicial.  Hence we can find $\tilde m
\in M$ such that
\d
\langle \tilde m,n_\rho\rangle > 0\ \hbox{when $\rho \in \Sigma(1) -
	A$, and\ } \langle \tilde m,n_\rho\rangle = 0\ \hbox{when
	$\rho \in A$}\ .
\d
Note that $\tilde m \in \check\sigma\cap\ker(i_A)$.  Now, given any $m
\in \ker(i_A)$, we can pick an integer $N$ such that $m + N\tilde m
\in \check\sigma\cap\ker(i_A)$.  Then $m = (m+N\tilde m) - N\tilde
m$ shows that $\check\sigma\cap\ker(i_A)$ generates $\ker(i_A)$ as a
group.  This proves that $g_A(\ker(i_A)) = \{1\}$.

	From the exact sequence (3), it follows that $g_A :
\z^{\Delta(1)-A} \to \cx^*$ induces a map ${\tilde g}_A : \tilde A \to
\cx^*$, where $\tilde A \subset A_{n-1}(X)$ is the image of
$\z^{\Delta(1)-A} \to A_{n-1}(X)$.  Since $\cx^*$ is injective as an
abelian group, ${\tilde g}_A$ extends to $g : A_{n-1}(X)
\to \cx^*$.  Thus $g \in G$, and it is straightforward to check that
$u_\rho = g([D_\rho])\,t_\rho$ for all $\rho$.  Then $u \in G\cdot t$,
and we conclude that $G\cdot t = V$.  Thus the orbits are closed and
we have a geometric quotient.

	To prove the converse, assume that $X$ is not simplicial.  The
following argument of I.~Musson will show that $G$ has a non-closed
orbit in $\cx^{\Delta(1)} - Z$.  Let $\sigma \in \Delta$ be a
non-simplicial cone, so that we have a relation $\sum_{\rho \in
\sigma(1)} b_\rho n_\rho = 0$.  We can assume that $b_\rho \in \z$ and
that at least one $b_\rho > 0$.  Then set $b_\rho = 0$ for $\rho \notin
\sigma(1)$.  Applying ${\rm Hom}_\z(-,\z)$ to the exact sequence at
the bottom of (1), we see that $\phi(t)([D_\rho]) = t^{b_\rho}$
defines an element of $G$ for $t \in \cx^*$ (in fact, $\phi(t)$ is a
$1$-parameter subgroup of $G$).  Now consider $u = (u_\rho) \in
\cx^{\Delta(1)}$, where
\d
u_\rho = \cases{1 & if $\rho \notin \sigma(1)$ or $m_\rho > 0$\cr
0 & otherwise\ .\cr}
\d
It is easy to check that $u$ and $\lim_{t \to 0}\phi(t)(u)$ are
in $\cx^{\Delta(1)}-Z$ but lie in different $G$-orbits.  Thus $G\cdot
u$ isn't closed in $\cx^{\Delta(1)}-Z$, and Theorem 2.1 is proved.\qed
\medskip

	The description of $X$ as a quotient of
$\cx^{\Delta(1)}-Z$ gives a new way of looking at the torus
action on $X$.  Namely, an element of the torus $(\cx^*)^{\Delta(1)}$
gives an automorphism of $\cx^{\Delta(1)}$ which preserves $Z$.
It also commutes with the action of $G$, so that we get an
automorphism of $X$ by the universal property of a categorical
quotient.  This gives a map $(\cx^*)^{\Delta(1)} \to {\rm Aut}(X)$,
and the kernel is easily seen to be $G$ (see \S4 for a proof).  Since
we have the exact sequence
\d
1 \longrightarrow G \longrightarrow (\cx^*)^{\Delta(1)}
	\longrightarrow T \longrightarrow 1\ ,
\d
it follows that the torus $T$ is naturally a subgroup of ${\rm
Aut}(X)$.  In \S4, we will discuss ${\rm Aut}(X)$ in more detail.
Note also that $(\cx^*)^{\Delta(1)} \subset \cx^{\Delta(1)} - Z$,
and thus, taking the quotient with respect to $G$, we get the
inclusion $T \subset X$.

	For the final result of this section, we will use Theorem 2.1
to study closed subsets (= reduced closed subschemes) of a simplicial
toric variety $X$.  Given a graded ideal $I \subset S$, note that the
set
\d
{\bf V}(I) - Z = \{ t \in \cx^{\Delta(1)} - Z : f(t) = 0\ \hbox{for
all $f \in I$}\} \subset \cx^{\Delta(1)}-Z
\d
is $G$-invariant and closed in the Zariski topology.  Since we have a
geometric quotient, the quotient map $\cx^{\Delta(1)}-Z \to X$ is
submersive (see [Fogarty and Mumford, Definition 0.6]).  Thus there is
a one-to-one correspondence between $G$-invariant closed sets in
$\cx^{\Delta(1)}-Z$ and closed sets in $X$.  Hence ${\bf V}(I) - Z$
determines a closed subset of $X$ which we will denote ${\bf V}_X(I)$.
By abuse of notation, we can write
\d
{\bf V}_X(I) = \{t\in X: f(t) =0\ \hbox{for all $f \in I$}\}\ ,
\d
provided we think of $t$ as the homogeneous coordinates of a point of
$X$.

	The map sending $I \subset S$ to ${\bf V}_X(I) \subset X$ has
the following properties:

\proclaim Proposition 2.4. Let $X$ be a simplicial toric variety, and
let $B = \langle x^{\hat\sigma}: \sigma \in \Delta\rangle \subset S$
be the ideal defined in \S1.
\vskip0pt
\item{(i)} (The Toric Nullstellensatz)  For any graded ideal $I
\subset S$, we have ${\bf V}_X(I) = \emptyset$ if and only if $B^m
\subset I$ for some integer $m$.
\item{(ii)} (The Toric Ideal--Variety Correspondence)  The map $I
\mapsto {\bf V}_X(I)$ induces a one-to-one correspondence
between radical graded ideals of $S$ contained in $B$ and closed
subsets of $X$.
\vskip0pt

\prf (i) Since ${\bf V}_X(I) = \emptyset$ in $X$ if and only if ${\bf
V}(I) \subset Z = {\bf V}(B)$ in $\cx^{\Delta(1)}$, the
first part of the proposition follows from the usual Nullstellensatz.

	(ii) Given a closed subset $W \subset \cx^{\Delta(1)}$, we get
the ideal ${\bf I}(W) = \{ f \in S: f_{|W} = 0\}$.  It is
straightforward to show that $W$ is $G$-invariant if and only if ${\bf
I}(W)$ is a graded ideal of $S$.  Note also that $B$ is radical since
it is generated by squarefree monomials.  Thus $B = {\bf I}(Z)$,
and we get one-to-one correspondences
\d
\eqalign{&\hbox{radical graded ideals of $S$ contained in $B$}\cr
	& \longleftrightarrow \hbox{$G$-invariant closed subsets of
	$\cx^{\Delta(1)}$ containing $Z$}\cr
	& \longleftrightarrow \hbox{$G$-invariant closed subsets of
	$\cx^{\Delta(1)}-Z$}\ ,\cr}
\d
where the last correspondence is induced by Zariski closure since $Z$
is $G$-invariant.  The proposition now follows since $\cx^{\Delta(1)}
-Z \to X$ is submersive.\qed
\medskip

	In \S3, we will generalize this proposition by showing that
there is a correspondence between graded $S$-modules and
quasi-coherent sheaves on $X$.
\bigskip

\noindent{\bf \S3. Graded Modules and Quasi-coherent Sheaves}
\bigskip

	In this section we will explore the relation between sheaves
on a toric variety $X$ and graded modules over its homogeneous
coordinate ring $S$.  The theory we will sketch is similar to what
happens for sheaves on $\P^n$ (as given in [Hartshorne] or [Serre]).

	An $S$-module $F$ is {\it graded\/} if there is a direct sum
decomposition
\d
F = \textstyle{\bigoplus_{\a \in A_{n-1}(X)}} F_\a
\d
such that $S_\a\cdot F_\beta \subset F_{\a+\beta}$ for all $\a,\beta
\in A_{n-1}(X)$.  Given such a module $F$, we will construct a
quasi-coherent sheaf $\widetilde F$ on $X$ as follows.  First, for
$\sigma \in \Delta$, let $S_\sigma$ be the localization of $S$ at
$x^{\hat\sigma} = \prod_{\rho \notin \sigma(1)} x_\rho$.  Then
$F_\sigma = F\otimes_S S_\sigma$ is a graded $S_\sigma$-module, so
that taking elements of degree $0$ gives a $(S_\sigma)_0$-module
$(F_\sigma)_0$.  From the proof of Theorem 2.1, we know that $X_\sigma
= {\rm Spec}((S_\sigma)_0)$ is an affine open subset of $X$, and then
standard theory tells us that $(F_\sigma)_0$ determines a
quasi-coherent sheaf $\widetilde{(F_\sigma)}_0$ on $X_\sigma$.

	To see how these sheaves patch together, let $\tau$ be a face
of $\sigma$.  Then $\tau = \sigma\cap\{m^\perp\}$ for some $m \in
\check\sigma\cap M$, and as in Lemma 2.3, one can prove
\d
(F_\sigma)_0 = \big((F_\tau)_0\big)_{x^{D_m}}
\d
where $x^{D_m} = \prod_\rho x_\rho^{\langle m,n_\rho\rangle}$.  It
follows immediately that the sheaves $\widetilde{(F_\sigma)}_0$ on
$X_\sigma$ patch to give a quasi-coherent sheaf $\widetilde{F}$ on
$X$.

	The map $F \mapsto \widetilde{F}$ has the following basic
properties of (we omit the simple proof):

\proclaim Proposition 3.1. The map sending $F$ to $\widetilde F$ is an
exact functor from graded $S$-modules to quasi-coherent ${\cal
O}_X$-modules.\qed

	As an example of how this works, consider the $S$-module
$S(\a)$, where $S(\a)_\beta = S_{\a+\beta}$ for all $\beta$.  We will
denote the sheaf $\widetilde{S(\a)}$ by ${\cal O}_X(\a)$.  If $\a =
[D]$ for some $T$-invariant divisor $D$, then there is an isomorphism
\d
{\cal O}_X(\a) \simeq {\cal O}_X(D)\ .\leqno{(5)}
\d
To prove this, it suffices to find compatible isomorphisms
\d
(S(\a)_\sigma)_0 \simeq H^0(X_\sigma,{\cal O}_X(D))
\d
for each affine open $X_\sigma$ of $X$.  If $D = \sum_\rho a_\rho
D_\rho$, then [Fulton, \S3.4] shows that this is equivalent to finding
compatible isomorphisms
\d
(S_\sigma)_\a \simeq \bigoplus_{\scriptstyle{{\scriptstyle{m \in M
\atop \langle m,n_\rho \rangle \ge -a_\rho}} \atop {\rm for\
all}\ \rho \in \sigma(1)}} \cx\cdot \chi^m\ .
\d
As in Proposition 1.1, we can map $\chi^m$ to the monomial $x^{D_m +
D} = \prod_\rho x^{\langle m,n_\rho\rangle + a_\rho} \in
(S_\sigma)_\a$ The compatibility of these maps is easily checked,
giving us the desired isomorphism (5).

	From (5), we see that $\widetilde{S} = {\cal O}_X$.  It
follows easily that $H^0(X,{\cal O}_X(\a))$ is canonically isomorphic
to $S_\a$, so that
\d
S \simeq \textstyle{\bigoplus_\a} H^0(X,{\cal O}_X(\a))\ .
\d
{}From (5), we also conclude that ${\cal O}_X(\a)$ is a coherent sheaf
(see [Fulton, \S3.4]).  However, it need not be a line bundle---in
fact, ${\cal O}_X(\a)$ is invertible if and only if $\a$ is the class
of a Cartier divisor.  See [Dolgachev, \S1.5] for an example of how
${\cal O}_X(\a)$ can fail to be invertible when $X$ is a weighted
projective space.  Also note that the natural map
\d
{\cal O}_X(\a)\otimes_{{\cal O}_X}\!{\cal O}_X(\beta) \to {\cal
O}_X(\a+\beta) \leqno(6)
\d
need not be an isomorphism---see [Dolgachev, \S1.5] for an example.
However, when ${\cal O}_X(\a)$ and ${\cal O}_X(\beta)$ are invertible,
then (6) is an isomorphism.

	Our next task is to see which sheaves on $X$ come from graded
$S$-modules.  When $X$ is simplicial, the answer is especially nice:

\proclaim Theorem 3.2.  If $X$ is a simplicial toric variety, then
every quasi-coherent sheaf on $X$ is of the form $\widetilde{F}$ for
some graded $S$-module $F$.

\prf The proof is similar to the proof of the corresponding result in
$\P^n$, as given in [Hartshorne] or [Serre].  The argument requires
standard facts about the behavior of sections of line bundles (see
[Hartshorne, Lemma II.5.14]).  As we observed above, ${\cal O}_X(\a)$
need not be a line bundle.  However, when $X$ is simplicial, ${\rm
Pic}(X)$ has finite index in $A_{n-1}(X)$ (this follows from an
exercise in [Fulton, \S3.3]).  Thus, given any $\a$, there is a
positive integer $k$ such that ${\cal O}_X(k\a)$ is a line bundle.

	Given a quasi-coherent sheaf ${\cal F}$ on $X$, let
\d
F = \textstyle{\bigoplus_\a} H^0(X,{\cal F}\otimes_{{\cal O}_X}\!{\cal
O}_X(\a))\ .
\d
Using (6) and the natural isomorphism $S_\a \simeq H^0(X,{\cal
O}_X(\a))$, we see that $F$ has the natural structure of a graded
$S$-module.  Given $\sigma \in \Delta$, let
\d
F'_\sigma = \textstyle{\bigoplus_\a} H^0(X_\sigma,{\cal
F}\otimes_{{\cal O}_X}\!{\cal O}_X(\a))
\d
which is a graded module over $S_\sigma$.  There is an obvious map
\d
\phi_\sigma : (F_\sigma)_0 \longrightarrow (F'_\sigma)_0 =
H^0(X_\sigma,{\cal F})\ .
\d
which patch to give $\phi : \widetilde{F} \to {\cal F}$.  To show that
$\phi$ is an isomorphism, it suffices to show that $\phi_\sigma$ is an
isomorphism for all $\sigma \in \Delta$.

	An element of $(F_\sigma)_0$ is written $f/x^D$, where $f \in
H^0(X,{\cal F}\otimes_{{\cal O}_X}\!{\cal O}_X(\a))$ and $\deg(x^D) =
\a$.  Since $f/x^D = (x^{(k-1)D}f)/x^{kD}$ when $k > 0$, we can assume
that ${\cal O}_X(\a)$ is a line bundle.  We can also assume that
$x^{\hat\sigma}$ divides $x^D$, which implies that $X_\sigma =
\{x \in X : x^D \ne 0\}$.  This guarantees that [Hartshorne, Lemma
II.5.14] applies, and from here, the proof follows the arguments in
[Hartshorne, Proposition II.5.15] without change. \qed
\medskip

	We next study coherent sheaves on $X$:

\proclaim Proposition 3.3.  Let $X$ be a toric variety with coordinate
ring $S$.
\vskip0pt
\item{(i)} If $F$ is a finitely generated graded $S$-module, the
$\widetilde{F}$ is a coherent sheaf on $X$.
\item{(ii)} If $X$ is simplicial, then every coherent sheaf on $X$ is
of the form $\widetilde{F}$ for some finitely generated graded
$S$-module $F$.
\vskip0pt

\prf (i) A set of homogeneous generators of $F$ gives a surjective map
$\oplus_i S(-\a_i) \to F$.  By Proposition 3.1, we get a surjection
$\oplus_i {\cal O}_X(-\a_i) \to \widetilde{F}$, and it follows that
$\widetilde{F}$ is coherent.

	(ii) Suppose that ${\cal F}$ is coherent.  On the affine
scheme $X_\sigma$, we can find finitely many sections $f_{i,\sigma}
\in H^0(X_\sigma,{\cal F})$ which generate ${\cal F}$ over $X_\sigma$.
Since $X$ is simplicial, arguing as in the proof of Theorem 3.2 shows
that we can find some $x^{D_\sigma}$ invertible in $S_\sigma$ such
that $x^{D_\sigma} f_{i,\sigma}$ comes from a section $g_{i,\sigma}$
of ${\cal F}\otimes_{{\cal O}_X}\!{\cal O}_X(\a_\sigma)$, where
$\a_\sigma = \deg(x^{D_\sigma})$.  Now consider the graded $S$-module
$F' \subset F = \oplus_\a H^0(X,{\cal F}\otimes_{{\cal O}_X}\!{\cal
O}_X(\a))$ generated by the $g_{i,\sigma}$.  Since $\widetilde{F}
\simeq {\cal F}$, we see that $\widetilde{F}' \to {\cal F}$ is
injective.  Over $X_\sigma$, we have $f_{i,\sigma} =
g_{i,\sigma}/x^{D_\sigma} \in (F'_\sigma)_0$, and since these sections
generate ${\cal F}$ over $X_\sigma$, it follows that $\widetilde{F}'
\simeq {\cal F}$.\qed
\medskip

	The map sending graded $S$-module $F$ to the sheaf
$\widetilde{F}$ is not injective, and in fact, there are many nonzero
modules which give the zero sheaf.  To analyze this phenomenon, we
will use the subgroup ${\rm Pic}(X) \subset A_{n-1}(X)$.  The
following lemma will be helpful:

\proclaim Lemma 3.4. If $\a \in {\rm Pic}(X) \subset A_{n-1}(X)$, then
for every $\sigma \in \Delta$, we can write $\a$ in the form $\a =
\bigl[\sum_{\rho \notin \sigma(1)} a_\rho D_\rho\bigr]$.

\prf It is well-known (see, for example, the exercises in [Fulton,
\S3.3]) that if $\sum_\rho b_\rho D_\rho$ is a Cartier divisor, then
for each $\sigma \in \Delta$, there is $m \in M$ such that $\langle
m,n_\rho\rangle = -b_\rho$ for all $\rho \in \sigma(1)$.  It follows
that $D$ is linearly equivalent to the divisor $\sum_{\rho
\notin \sigma(1)} (b_\rho +\langle m,n_\rho\rangle)D_\rho$, and the
lemma is proved.\qed
\medskip

	We can now describe which finitely generated graded
$S$-modules $F$ correspond to the zero sheaf on $X$.  In the
projective case, where $S = \cx[x_0,\dots,x_n]$, this happens if and
only if $F_k = \{0\}$ for all $k \gg 0$.  This can be rephrased in
terms of the irrelevant ideal $\langle x_0,\dots,x_n\rangle \subset
S$: namely, $F_k = \{0\}$ for $k \gg 0$ exactly when $\langle
x_0,\dots,x_n\rangle^k F = \{0\}$ for some $k$.  In the case of a
simplicial toric variety, we can generalize this as follows:

\proclaim Proposition 3.5. Let $X$ be a simplicial toric variety, and
let $B \subset S$ be the ideal defined in \S1.  Then, for a finitely
generated graded $S$-module $F$, we have $\widetilde{F} = 0$ on $X$ if
and only if there is some $k > 0$ such that $B^k F_\a = \{0\}$ for all
$\a \in {\rm Pic}(X)$.

\prf First, suppose that $B^k F_\a = \{0\}$ for all $\a \in {\rm
Pic}(X)$.  Then fix $\sigma \in \Delta$ and take $f/x^D \in
(F_\sigma)_0$.  Because $X$ is simplicial, we know that $l[D] \in {\rm
Pic}(X)$ for some $l > 0$.  If we let $\a = l[D]$, then $x^{(l-1)D}f
\in F_\a$, and hence $(x^{\hat\sigma})^k x^{(l-1)D}f = 0$ for some $k
> 0$.  This easily implies $f/x^D = 0$, and it follows that
$(F_\sigma)_0 = \{0\}$.  Thus $\widetilde{F}$ is the zero sheaf.

	Conversely, let $F$ be a finitely generated graded $S$-module
such that $\widetilde{F} = 0$.  Fix $\a \in {\rm Pic}(X)$ and let $f
\in F_\a$.  By Lemma 3.4, we can write $\a = [D]$ where $D = \sum_{\rho
\notin \sigma(1)} a_\rho D_\rho$.  Then $x^D$ is an invertible element
in $S_\sigma$, and hence $f/x^D \in (F_\sigma)_0 = \{0\}$.  This
implies that $(x^{\hat\sigma})^k f = 0$ for some $k > 0$.  To show
that we can get a $k$ that works for all $f \in F_\a$ and all $\a \in
{\rm Pic}(X)$, note that $\sum_{\a \in {\rm Pic}(X)} S F_\a \subset F$
is a submodule and hence is finitely generated.  The proposition now
follows easily.\qed
\medskip

	When $X$ is a smooth toric variety, we know that ${\rm Pic}(X)
= A_{n-1}(X)$.  Hence we get the following corollary of Proposition
3.5:

\proclaim Corollary 3.6. If $X$ is a smooth toric variety and $F$ is
a finitely generated graded $S$-module, then $\widetilde{F} = 0$ on
$X$ if and only if there is some $k > 0$ with $B^k F = \{0\}$.\qed

	To see how Corollary 3.6 can fail when $X$ not smooth,
consider $X = \P(1,1,2)$.  Here, $S = \cx[x,y,z]$, where $x,y$ have
degree $1$ and $z$ has degree $2$, and $B = \langle x,y,z\rangle$.
The graded $F = S(1)/(xS(1)+yS(1))$ has only elements of odd degree.
Then $(F_z)_0 = \{0\}$ since $z$ has degree $2$, and it is clear that
$(F_x)_0 = (F_y)_0 = \{0\}$.  It follows that $\widetilde{F} = 0$, yet
one easily checks that $B^k F = z^k F \ne \{0\}$ for all $k$.  Thus
Corollary 3.6 does not hold.

	We will end this section with a study of the ideal-variety
correspondence.  Our goal is to refine what we did in Proposition 2.4
by taking the scheme structure into account.  First note that if $I
\subset S$ is a graded ideal, then by exactness we get an ideal sheaf
$\widetilde{I} \subset {\cal O}_X$.  Hence $I \subset S$ determines
the closed subscheme $Y \subset X$ with ${\cal O}_Y = {\cal
O}_X/\widetilde{I}$.  The following theorem tells us which subschemes
of $X$ arise in this way:

\proclaim Theorem 3.7. Let $X$ be a simplicial toric variety, and let
$B \subset S$ be as usual.
\vskip0pt
\item{(i)} Every closed subscheme of $X$ comes from a graded ideal of
$S$ (as described above).
\item{(ii)} Two graded ideals $I$ and $J$ of $S$ correspond to the
same closed subscheme of $X$ if and only $(I\colon\! B^\infty)_\a =
(J\colon\! B^\infty)_\a$ for all $\a \in {\rm Pic}(X)$.
\vskip0pt

\prf (i) Given an ideal sheaf ${\cal I} \subset {\cal O}_X$, let $I' =
\oplus_\a H^0(X,{\cal I}\otimes_{{\cal O}_X}\!{\cal O}_X(\a))$.  Since
$X$ is simplicial, we know that $\widetilde{I}' = {\cal I}$.  The map
$I' \to S$ factors $I' \ontoarrow I \hookrightarrow S$, which gives
maps ${\cal I} \ontoarrow \widetilde{I} \hookrightarrow {\cal O}_X$ by
exactness.  Since the composition is ${\cal I} \subset {\cal O}_X$, it
follows that ${\cal I} = \widetilde{I}$.

	(ii) First, note that if $I$ is graded, then so is $I\colon\!
B^\infty$.  Now, given graded ideals $I,J \subset S$, assume that
$\widetilde{I} = \widetilde{J}$.  Then for every $\sigma \in \Delta$,
we have $(I_\sigma)_0 = (J_\sigma)_0$ in $(S_\sigma)_0$.  Let $\a \in
{\rm Pic}(X)$ and pick $u \in (I\colon\!  B^\infty)_\a$.  For $k$
sufficiently large, we have $(x^{\hat\sigma})^k u \in I$.  Increasing
$k$ if necessary, we can assume that $\beta = k\deg(x^{\hat\sigma})
\in {\rm Pic}(X)$.  Thus $\a + \beta \in {\rm Pic}(X)$, so Lemma 3.4
implies $\a + \beta = [D]$, where $D = \sum_{\rho \notin \sigma(1)}
a_\rho D_\rho$.  Then $((x^{\hat\sigma})^k u)/x^D \in (I_\sigma)_0 =
(J_\sigma)_0$, which implies that $(x^{\hat\sigma})^{k+l} u \in J$ for
some $l > 0$.  This is true for all $\sigma$, so that $u \in
(J\colon\!  B^\infty)_\a$.  The other inclusion $(J\colon\!
B^\infty)_\a \subset (I\colon\!  B^\infty)_\a$ follows by symmetry,
and equality is proved.

	Conversely, suppose $(I\colon\! B^\infty)_\a = (J\colon\!
B^\infty)_\a$ for all $\a \in {\rm Pic}(X)$.  To prove $\widetilde{I}
= \widetilde{J}$, it suffices to show that $(I_\sigma)_0 =
(J_\sigma)_0$ for all $\sigma$.  But given $u \in (I_\sigma)_0$, we
know that $u = f/(x^{\hat\sigma})^k$, and we can assume that $\a =
\deg(f) \in {\rm Pic}(X)$.  Then $f \in I_\a \subset (I\colon\!
B^\infty)_\a = (J\colon\!  B^\infty)_\a$.  This easily implies $u =
f/(x^{\hat\sigma})^k$ is in $(J_\sigma)_0$.  Thus $(I_\sigma)_0
\subset (J_\sigma)_0$, and the opposite inclusion follows by symmetry.
This completes the proof of the theorem. \qed
\medskip

	As an example, consider the ideal $I = \langle x_\rho\rangle$
for some fixed $\rho \in \Delta(1)$.  We will show that the
corresponding subscheme of $X$ is the divisor $D_\rho$.  Given $\sigma
\in \Delta$, we need to study $(I_\sigma)_0$.  If $\rho
\in \sigma(1)$, then it is easy to check that under the isomorphism
$(S_\sigma)_0 \simeq \cx[\check\sigma\cap M]$, we have
\d
(I_\sigma)_0 \simeq
\displaystyle{\bigoplus_{\scriptstyle{{\scriptstyle{m \in M
\atop \langle m,\sigma \rangle \ge 0}} \atop {\langle m,n_\rho
\rangle > 0}}}} \cx\cdot \chi^m \ .
\d
By an exercise in [Fulton, \S3.1], this is exactly the ideal of
$D_\rho$ in the affine piece $X_\sigma \subset X$.  On the other hand,
if $\rho \notin \sigma(1)$, then $(I_\sigma)_0 = \langle 1\rangle$,
which is consistent with $D_\rho \cap X_\sigma = \emptyset$.  Thus
$D_\rho$ is defined scheme-theoretically by the ideal $\langle
x_\rho\rangle$.  In particular, $D_\rho$ has the global equation
$x_\rho = 0$ even though it may not be a Cartier divisor.

	In the smooth case, Theorem 3.7 has some nice consequences.
We say that an ideal $I \subset S$ is {\it $B$-saturated\/} if
$I\colon\!  B = I$.  Then we have:

\proclaim Corollary 3.8. Let $X$ be a smooth toric variety.
\vskip0pt
\item{(i)} Two graded ideals $I$ and $J$ of $S$ give the same closed
subscheme of $X$ if and only if $I\colon\! B^\infty = J\colon\!
B^\infty$.
\item{(ii)} There is a one-to-one correspondence between closed
subschemes of $X$ and graded $B$-saturated ideals of $S$.

\prf Part (i) follows from Theorem 3.7 since ${\rm Pic}(X) =
A_{n-1}(X)$, and part (ii) follows immediately since an ideal is
$B$-saturated if and only if it is of the form $J\colon\! B^\infty$.
\qed
\medskip

	When $X$ is smooth, Corollary 3.8 gives a one-to-one
correspondence between radical graded $B$-saturated ideals of $S$ and
closed reduced subschemes of $X$.  This is slightly different from the
correspondence given in Proposition 2.4, where we used radical ideals
contained in $B$.  However, the maps
\d
\eqalign{\hbox{$I$ radical, graded, $B$-saturated}\ &\mapsto\  I\cap
B\cr
\hbox{$J$ radical, graded, contained in $B$}\ &\mapsto\ J\colon\!B\cr}
\d
give bijections between these two classes of radical graded ideals, so
the two descriptions of reduced subschemes are compatible.

	Finally, one could ask if there is an analog of Corollary 3.8
in the simplicial case.  There are two ways to approach this problem.
One way is to work with the following special class of graded ideals.
We say that a graded ideal $I \subset S$ is {\it Pic-generated\/} if
$I$ is generated by homogeneous elements whose degrees lie in ${\rm
Pic}(X)$, and $I$ is {\it Pic-saturated\/} if $I_\a = (I\colon\!
B)_\a$ for all $\a \in {\rm Pic}(X)$.  For example, if $I$ is any
graded ideal, it is easy to check that $\sum_{\a \in {\rm
Pic}(X)}S(I\colon\!  B^\infty)_\a$ is both Pic-generated and
Pic-saturated.  Then one can prove the following result (we will omit
the proof):

\proclaim Corollary 3.9.  If $X$ is a simplicial toric variety, then
there is a one-to-one correspondence between closed subschemes of $X$
and graded ideals of $S$ which are Pic-generated and
Pic-saturated.\qed

	Another way to study sheaves on $X$ is to replace the ring $S$
with the subring
\d
R = \textstyle{\bigoplus_{\a \in {\rm Pic}(X)}} S_\a\ .\leqno{(7)}
\d
When $X$ is simplicial, we know that ${\rm Pic}(X) \subset
A_{n-1}(X)$ has finite index.  Thus
\d
H = {\rm Hom}_\z(A_{n-1}(X)/{\rm Pic}(X),\cx^*) \subset G = {\rm
Hom}_\z(A_{n-1}(X),\cx^*)
\d
is a finite subgroup of $G$, and $R = S^H$ under the action of $H$ on
$S$.  It follows that $Y = {\rm Spec}(R)$ is the quotient of
$\cx^{\Delta(1)}$ under the action of $H$.  Furthermore, if we set $B'
= B\cap R$ and let $Z' = {\bf V}(B') \subset Y$, then it is easy to
see that $Y - Z' = (\cx^{\Delta(1)} - Z)/H$.  There is also a natural
action of ${\rm Hom}_\z({\rm Pic}(X),\cx^*)$ on $Y-Z'$, and Theorem
2.1 implies the following result:

\proclaim Theorem 3.10. If $X$ is a simplicial toric variety and $Y =
{\rm Spec}(R)$, where $R$ is the ring (7), then $X$ is the geometric
quotient of the action of ${\rm Hom}_\z({\rm Pic}(X),\cx^*)$ on
$Y-Z'$.\qed

	Although the ring $R$ may not be a polynomial ring, the
relation between sheaves and $R$-modules is nicer than the relation
between sheaves and $S$-modules.  It is easy to show that a graded
$R$-module $F$ (where the grading is now given by ${\rm Pic}(X)$)
gives a sheaf $\widetilde{F}$ on $X$.  Then one can prove the
following result (we omit the proof since the arguments are similar to
what we did above):

\proclaim Theorem 3.11.  If $X$ is a simplicial toric variety and $R$
is the ring defined in (7), then:
\vskip0pt
\item{(i)} The map sending a graded $R$-module $F$ to $\widetilde{F}$
is an exact functor.
\item{(ii)} All quasi-coherent sheaves on $X$ arise in this way.
\item{(iii)} $\widetilde{F}$ is coherent when $F$ is finitely
generated, and this gives all coherent sheaves on $X$.
\item{(iv)} If $F$ is finitely generated, then $\widetilde{F} = 0$ if
and only if ${B'}^k F = \{0\}$ for some $k > 0$ (where $B' = B\cap
R$).\qed
\vskip0pt

	Finally, a graded ideal $I \subset R$ gives an ideal sheaf
$\widetilde{I} \subset {\cal O}_X$, so that $I$ determines a closed
subscheme of $X$.  We will omit the proof of the following result:

\proclaim Theorem 3.12. If $X$ is a simplicial toric variety and $R$
is the ring defined in (7), then:
\vskip0pt
\item{(i)} Every closed subscheme of $X$ comes from some graded ideal
of $R$.
\item{(ii)} Two graded ideals $I,J$ give the same closed subscheme if
and only if $I\colon\! {B'}^\infty = J\colon\! {B'}^\infty$ (where $B'
= B\cap R$).
\item{(iii)} There is a one-to-one correspondence between graded
$B'$-saturated ideals of $R$ and closed subschemes of $X$.\qed
\vskip0pt

\bigskip

\noindent{\bf \S4. The Automorphism Group}
\bigskip

	This section will use the graded ring $S$ to give an explicit
description of the automorphism group ${\rm Aut}(X)$ of a complete
simplicial toric variety $X$.  To see what sort of results to expect,
consider $\P^n \simeq (\cx^{n+1}-\{0\})/\cx^*$.  Notice that ${\rm
Aut}(\P^n) = PGL(n+1,\cx)$ doesn't act on $\cx^{n+1}-\{0\}$; rather,
we must use $GL(n+1,\cx)$, which is related to $PGL(n+1,\cx)$ by the
exact sequence
\d
1 \longrightarrow \cx^* \longrightarrow GL(n+1,\cx)
\longrightarrow PGL(n+1,\cx) \longrightarrow 1\ .
\d
For a complete simplicial toric variety $X$, we know $X \simeq
(\cx^{\Delta(1)}-Z)/G$ by Theorem~2.1.  Then the analog of
$GL(n+1,\cx)$ will be the following pair of groups:

\proclaim Definition 4.1.
\vskip0pt
\item{(i)} $\widetilde{\rm Aut}(X)$ is the normalizer
of $G$ in the automorphism group of $\cx^{\Delta(1)} -Z$.
\item{(ii)} $\widetilde{\rm Aut}{}^0(X)$ is the centralizer
of $G$ in the automorphism group of $\cx^{\Delta(1)} -Z$.
\vskip0pt

It is clear that $\widetilde{\rm Aut}{}^0(X) \subset \widetilde{\rm
Aut}(X)$ is a normal subgroup.  Note also that every element of
$\widetilde{\rm Aut}(X)$ preserves $G$-orbits and hence descends to an
automorphism of $X$.  Thus we get a natural homomorphism
$\widetilde{\rm Aut}(X) \to {\rm Aut}(X)$.

	We can now state the main result of this section:

\proclaim Theorem 4.2. Let $X$ be a complete simplicial toric variety,
and let $S = \cx[x_\rho]$ be its homogeneous coordinate ring.  Then:
\vskip0pt
\item{(i)} $\widetilde{\rm Aut}(X)$ is an affine algebraic group of
dimension $\sum_\rho \dim_\cx S_{\deg(x_\rho)}$, and
$\widetilde{\rm Aut}{}^0(X)$ is the connected component of the identity.
\item{(ii)} The natural map $\widetilde{\rm Aut}(X) \to {\rm Aut}(X)$
induces an exact sequence
\d
1 \longrightarrow G \longrightarrow \widetilde{\rm Aut}(X)
\longrightarrow {\rm Aut}(X) \longrightarrow 1\ .
\d
\item{(iii)} $\widetilde{\rm Aut}{}^0(X)$ is naturally isomorphic to the
group ${\rm Aut}_g(S)$ of graded $\cx$-algebra automorphisms of
$S$.
\vskip0pt

	In the smooth case, this theorem follows from results of
[Demazure].  Before we can begin the proof of Theorem 4.2, we need to
study the group ${\rm Aut}_g(S)$.  Consider the set ${\rm End}_g(S)$
of graded $\cx$-algebra homomorphisms $\phi:S\to S$ which satisfy
$\phi(1) = 1$.  Since $X$ is complete, we know that $S_0 = \cx$, and
it follows that $S^+ = \sum_{\a \ne 0} S_\a$ is closed under
multiplication.  Thus there is a natural bijection ${\rm End}_g(S)
\simeq {\rm End}_g(S^+)$.  Since $S^+$ has no identity to worry about,
${\rm End}_g(S^+)$ has a natural structure as a $\cx$-algebra, and in
this way we can regard ${\rm End}_g(S)$ as a $\cx$-algebra.

	Note that an element $\phi \in {\rm End}_g(S)$ is uniquely
determined by $\phi(x_\rho) \in S_{\deg(x_\rho)}$ for $\rho \in
\Delta(1)$.  Furthermore, each $S_{\deg(x_\rho)}$ is finite
dimensional since $X$ is complete, so that ${\rm End}_g(S)$ is a
finite dimensional $\cx$-algebra.  Since ${\rm Aut}_g(S)$ is the set
of invertible elements of ${\rm End}_g(S)$, standard arguments show
that ${\rm Aut}_g(S)$ is an affine algebraic group.

	To describe the structure of this group in more detail, we
will use the equivalence relation on $\Delta(1)$ defined by
\d
\rho \sim \rho' \iff \deg(x_\rho) = \deg(x_{\rho'})\ {\rm in}
\ A_{n-1}(X)\ .
\d
This partitions $\Delta(1)$ into disjoint subsets $\Delta_1\cup\cdots
\cup\Delta_s$, where each $\Delta_i$ corresponds to a set of variables
of the same degree $\a_i$.  Then each $S_{\a_i}$ can be written as
\d
S_{\a_i} = S_{\a_i}'\oplus S_{\a_i}''\leqno{(7)}
\d
where $S_{\a_i}'$ is spanned by the $x_\rho$ for $\rho \in \Delta_i$,
and $S_{\a_i}''$ is spanned by the remaining monomials in $S_{\a_i}$
(all of which are the product of at least two variables).

	Before we state our next result, note that the torus
$(\cx^*)^{\Delta(1)}$ can be regarded as a subgroup of ${\rm
Aut}_g(S)$ since $t = (t_\rho) \in (\cx^*)^{\Delta(1)}$ gives the
graded automorphism $\phi_t :S\to S$ defined by $\phi_t(x_\rho) =
t_\rho x_\rho$.  Then the group structure of ${\rm Aut}_g(S)$ can be
described as follows:

\proclaim Proposition 4.3.  Let $X$ be a complete toric variety, and
let $S$ be its homogeneous coordinate ring.  Then:
\vskip0pt
\item{(i)} ${\rm Aut}_g(S)$ is a connected affine algebraic group
of dimension $\sum_{i=1}^s |\Delta_i|\dim_\cx S_{\a_i}$, and
$(\cx^*)^{\Delta(1)} \subset {\rm Aut}_g(S)$ is its maximal torus.
\item{(ii)} The unipotent radical $R_u$ of ${\rm Aut}_g(S)$ is
isomorphic as a variety to an affine space of dimension $\sum_{i=1}^s
|\Delta_i|(\dim_\cx S_{\a_i} - |\Delta_i|)$.
\item{(iii)} ${\rm Aut}_g(S)$ contains a subgroup $G_s$ isomorphic
to the reductive group $\prod_{i=1}^s GL(S_{\a_i}')$ of dimension
$\sum_{i=1}^s |\Delta_i|^2$.
\item{(iv)} ${\rm Aut}_g(S)$ is isomorphic to the semidirect
product $R_u \semidir G_s$.
\vskip0pt

\prf We first study the structure of the $\cx$-algebra ${\rm
End}_g(S)$.  We have a natural vector space isomorphism
\d
{\rm End}_g(S) \simeq \oplus_{i=1}^s {\rm
Hom}_\cx(S_{\a_i}',S_{\a_i})
\d
since $S_{\a_i}'$ is spanned by $x_\rho$ for $\rho \in \Delta_i$ and a
graded homomorphism $\phi : S \to S$ is uniquely determined by
$\phi(x_\rho) \in S_{\a_i}$.  Then the direct sum decomposition (7)
gives an exact sequence
\d
0 \longrightarrow \oplus_{i=1}^s {\rm Hom}_\cx(S_{\a_i}',S_{\a_i}'')
\longrightarrow {\rm End}_g(S)\ \mapname{\Psi}\, \oplus_{i=1}^s {\rm
End}_\cx(S_{\a_i}') \longrightarrow 0\ .\leqno{(8)}
\d
Notice that $\Psi$ is a ring homomorphism.  This follows immediately
from the observation that $\phi(S_{\a_i}'') \subset S_{\a_i}''$ for
any $\phi \in {\rm End}_g(S)$.  Thus we can regard ${\cal N} =
\oplus_{i=1}^s {\rm Hom}_\cx (S_{\a_i}',S_{\a_i}'')$ as an ideal of
${\rm End}_g(S)$.  Note that $\phi \in {\cal N}$ if and only if
$\phi(x_\rho) \in S_{\a_i}''$ for $\rho \in \Delta_i$.

	We claim that every element of ${\cal N}$ is nilpotent.  To
prove this, we will use the order relation on the monomials of $S$
defined in Lemma 1.3.  For variables, the relation is defined by
$x_{\rho'} < x_{\rho}$ if and only if there is a monomial $x^D$ with
$\deg(x^D) = \deg(x_\rho)$, $x_{\rho'} | x^D$, and $x_{\rho'} \ne
x^D$.  If $\deg(x_\rho) = \a_i$, this is equivalent to the existence
of $x^D \in S_{\a_i}''$ with $x_{\rho'} | x^D$.

	Given $\phi \in {\cal N}$, let $N(\phi)$ the be number of
variables $x_\rho$ for which $\phi(x_\rho) = 0$.  If $N(\phi) =
|\Delta(1)|$, then $\phi = 0$ and we are done.  If $N(\phi) <
|\Delta(1)|$, then look at the variables where $\phi$ doesn't vanish
and pick one, say $x_\rho$, which is {\it minimal\/} with respect to
$<$.  We can do this because $X$ is complete and hence $<$ is
antisymmetric by Lemma 1.3.  If $\deg(x_\rho) = \a_i$, then as we
noted above, $x_{\rho'} < x_\rho$ for any variable $x_{\rho'}$
dividing a monomial in $S_{\a_i}''$.  By the minimality of $x_\rho$,
it follows that $\phi(x_{\rho'}) = 0$, and consequently
$\phi(S_{\a_i}'') = 0$.  This implies that $\phi^2(x_\rho) = 0$ since
$\phi(x_\rho) \in S_{\a_i}''$.  We conclude that $N(\phi^2) > N(\phi)$
since $\phi^2$ also vanishes on every variable killed by $\phi$.  The
same argument applies to $\phi^2, \phi^4, \phi^8, \dots \in {\cal N}$,
and it follows that some power of $\phi$ must be zero.  Thus $\phi$ is
nilpotent as claimed.

	Once we know that ${\cal N}$ is a nilpotent ideal, it follows
that every element of $1 + {\cal N}$ is invertible, so that $1+{\cal
N}$ is a unipotent subgroup of ${\rm Aut}_g(S)$.  From the exact
sequence (8), we also see that a element $\phi \in {\rm End}_g(S)$
is invertible if and only if its image in $\prod_{i=1}^s
GL(S_{\a_i}')$ is invertible.  From this, it follows that we get an
exact sequence
\d
1 \longrightarrow 1 + {\cal N}
\longrightarrow {\rm Aut}_g(S) \longrightarrow
\textstyle{\prod_{i=1}^s} GL(S_{\a_i}') \longrightarrow 1.\leqno(9)
\d
Since $\prod_{i=1}^s GL(S_{\a_i}')$ is a reductive group, it follows
that $1+{\cal N}$ is the unipotent radical of ${\rm Aut}_g(S)$.
Then part (ii) of the proposition follows because as a variety, we have
$1+{\cal N} \simeq {\cal N}$, which is an affine space of the required
dimension.

	To prove part (iii), observe that we have a map
\d
s: \textstyle{\prod_{i=1}^s} {\rm End}_\cx(S_{\a_i}') \to {\rm
End}_g(S)
\d
where $(\phi_i) \in \prod_{i=1}^s {\rm End}_\cx (S_{\a_i}')$ maps to
the homomorphism $\phi$ defined by $\phi(x_\rho) = \phi_i(x_\rho)$
when $\rho \in \Delta_i$.  It is easy to check that $s$ is a ring
homomorphism which is a section of the map $\Psi$ of (8).  It follows
that $s$ induces a map
\d
s^*:\textstyle{\prod_{i=1}^s} GL(S_{\a_i}') \longrightarrow {\rm
Aut}_g(S)
\d
which is a section of the exact sequence (9).  The image $G_s$ of
$s^*$ is a subgroup of ${\rm Aut}_g(S)$ isomorphic to $\prod_{i=1}^s
GL(S_{\a_i}')$, and thus ${\rm Aut}_g(S)$ is the semidirect product of
$(1+{\cal N})\semidir G_s$.  This proves parts (iii) and (iv), and
then part (i) now follows immediately.  \qed
\medskip

	We will next study the {\it roots\/} of ${\rm Aut}(X)$ as
defined by Demazure.  Recall that the roots of the complete toric
variety $X$ are the set
\d
R(N,\Delta)=\{m \in M: \hbox{$\exists \rho \in \Delta(1)$ with
$\langle m,n_\rho\rangle = 1$ and $\langle m,n_{\rho'}\rangle \le 0$
for $\rho' \ne \rho$}\}
\d
(see [Demazure, \S3.1] or [Oda, \S3.4]).  We divide the roots
$R(N,\Delta)$ into two classes as follows:
\d
R_s(N,\Delta) = R(N,\Delta)\cap(-R(N,\Delta))\, \quad R_u(N,\Delta) =
R(N,\Delta) - R_s(N,\Delta)\ .
\d
The roots in $R_s(N,\Delta)$ are called the semisimple roots.

	We can describe $R(N,\Delta)$ and $R_s(N,\Delta)$ in terms of
the graded ring $S$ as follows:

\proclaim Lemma 4.4. There is are one-to-one correspondences
\d
\eqalign{R(N,\Delta) &\longleftrightarrow \{(x_\rho,x^D) : \hbox{$x^D
\in S$ is a
monomial, $x^D \ne x_\rho$, $\deg(x^D) = \deg(x_\rho)$}\}\cr
R_s(N,\Delta) &\longleftrightarrow \{(x_\rho,x_{\rho'}) : x_\rho \ne
x_{\rho'},\ \deg(x_\rho) = \deg(x_{\rho'})\}\ .\cr}
\d

\prf Given $m \in R(N,\Delta)$, we have $\rho$ with $\langle
m,n_\rho\rangle = 1$ and $\langle m,n_{\rho'}\rangle \le 0$ for $\rho'
\ne \rho$.  Then consider the map
\d
m \mapsto \biggl(x_\rho,\prod_{\rho' \ne \rho} x_{\rho'}^{\langle
-m,n_{\rho'}\rangle}\biggr)\ .
\d
Since $m$ is a root, the second component has nonnegative exponents
and hence is a monomial $x^D \in S$.  It is also clear that $x_\rho$
and $x^D$ are distinct and have the same degree.

	This map is injective because $M \to \z^{\Delta(1)}$ is
injective.  To see that it is surjective, suppose we have $x_\rho \ne
x^D$ of the same degree in $S$.  A first observation is that $x_\rho
\notdiv x^D$, for otherwise we could write $x^D = x^\rho x^E$, and it
would follow that $x^E \in S$ would have degree zero.  Yet $X$
complete implies $S_0 = \cx$, so that $x^E$ must be constant.  This
would force $x_\rho = x^D$, which is impossible.  Then, writing $x^D =
\prod_{\rho' \ne \rho} x_{\rho'}^{a_{\rho'}}$, the definition of
$\deg(x_\rho) = \deg(x^D)$ implies there is $m \in M$ with
\d
1 = \langle m,n_\rho\rangle + 0\ ,\quad 0 = \langle
m,n_{\rho'}\rangle + a_{\rho'}\ \hbox{for $\rho' \ne \rho$}\ .
\d
Thus $m$ is a root mapping to the pair $(x_\rho,x^D)$.

	The correspondence for semisimple roots follows easily from
the observation that if $m \leftrightarrow (x_\rho,x^D)$ and $-m
\leftrightarrow (x_{\rho'}, x^{D'})$, then $x^D = x_{\rho'}$. \qed
\medskip

	This lemma allows us to think of semisimple roots as pairs of
distinct variables of the same degree.

	Given a root $m \in R(N,\Delta)$, we have the corresponding
pair $(x_\rho,x^D)$.  Then, for $\lambda \in \cx$, consider the map
\d
y_m(\lambda) : S \longrightarrow S
\d
defined by $y_m(\lambda)(x_\rho) = x_\rho + \lambda x^D$ and
$y_m(\lambda)(x_{\rho'}) = x_{\rho'}$ for $\rho' \ne \rho$.  It is
easy to see that $y_m(\lambda)$ gives a one parameter group of graded
automorphisms of $S$, so that $y_m(\lambda) \in {\rm Aut}_g(S)$.  We
can relate the $y_m(\lambda)$ to ${\rm Aut}_g(S)$ as follows:

\proclaim Proposition 4.5. Let $X$ be a complete toric
variety.  Then ${\rm Aut}_g(S)$ is generated by the torus
$(\cx^*)^{\Delta(1)}$ and the one parameter subgroups $y_m(\lambda)$
for $m \in R(N,\Delta)$.

\prf In Proposition 4.3, we gave showed that ${\rm Aut}_g(S)$ was
the semidirect product
\d
{\rm Aut}_g(S) \simeq (1+{\cal N}) \semidir G_s\ ,\ {\rm where}\
G_s \simeq \textstyle{\oplus_{i=1}^s GL(S_{\a_i}')}\ .
\d
In this decomposition, we know from Proposition 4.3 that
$(\cx^*)^{\Delta(1)} \subset {\rm Aut}_g(S)$ is the maximal torus and
hence lies in $G_s$.

	First suppose that $m \leftrightarrow (x_\rho,x_{\rho'})$ is
in $R_s(N,\Delta)$.  Then, letting $\a_i = \deg(x_\rho)$, we have
$x_\rho, x_{\rho'} \in S_{\a_i}'$ in the direct sum (7).  Since
$y_m(\lambda)(x_\rho) = x_\rho + \lambda x_{\rho'}$ and
$y_m(\lambda)(x_{\rho''}) = x_{\rho''}$ for $\rho'' \ne \rho$, it
follows that $y_m(\lambda) \in G_s$, and in fact, the $y_m(\lambda)$
correspond to the roots of $G_s = \oplus_{i=1}^s GL(S_{\a_i}')$ in the
usual sense of algebraic groups.  Since any reductive group is
generated by its maximal torus and the one parameter subgroups given
by its roots, it follows that $G_s$ is generated by
$(\cx^*)^{\Delta(1)}$ and the $y_m(\lambda)$ for $m \in
R_s(N,\Delta)$.

	Next, suppose that $m \leftrightarrow (x_\rho, x^D)$ is a root
in $R_u(N,\Delta)$.  If we set $\a_i = \deg(x_\rho)$, then $x^D \in
S_{\a_i}''$ in the direct sum (7).  It follows that $y_m(\lambda) \in
1+{\cal N}$.  We claim that $1+{\cal N}$ is generated by the
$y_m(\lambda)$ for $m \in R_u(N,\Delta)$.  To prove this, suppose that
$\phi \in 1+{\cal N}$.  If $\phi(x_\rho) = x_\rho + u_\rho$, where
$u_\rho \in S_{\a_i}''$, then define the {\it weight\/} of $\phi$ to
be the sum
\d
w(\phi) = \sum_\rho \hbox{number of nonzero coefficients of $u_\rho$}\ .
\d
Note the $\phi$ is the identity if and only if $w(\phi) = 0$.  If
$\phi$ is not the identity, then we can use the ordering of variables
from the proof of Proposition 4.3 to pick a $x_\rho$ with $u_\rho \ne
0$ such that $x_\rho$ is {\it maximal\/} with respect to this property.
Let $\lambda x^D$ is a nonzero term of $u_\rho$, so that $u_\rho =
\lambda x^D + u'_\rho$.  Then $(x_\rho,x^D)$ corresponds to a root $m
\in R_u(N,\Delta)$, and we have
\d
y_m(-\lambda)\circ \phi(x_\rho) = y_m(-\lambda)(x_\rho + \lambda x^D +
u'_\rho) = x_\rho + u'_\rho \leqno{(10)}
\d
since none of the monomials in $u_\rho = \lambda x^D + u'_\rho$
involve $x_\rho$.  Next suppose that $\rho' \ne \rho$.  Then $x_\rho$
doesn't appear in any monomial in $u_{\rho'}$, for if it did, we would
have $x_\rho < x_{\rho'}$, which would contradict the maximality of
$x_\rho$.  Thus $\phi(x_{\rho'}) = x_{\rho'} + u_{\rho'}$ only
involves variables different from $x_\rho$, and hence
\d
y_m(-\lambda)\circ\phi(x_{\rho'}) = \phi(x_{\rho'})\quad \hbox{for
$\rho' \ne \rho$}\ .\leqno{(11)}
\d
{}From (10) and (11) we conclude that $w(y_m(-\lambda)\circ\phi) <
w(\phi)$, and from here it is easy to see $\phi$ can be written in
terms of the $y_m(\lambda)$.  This shows that $1+{\cal N}$ is
generated by the $y_m(\lambda)$ for $m \in R_u(N,\Delta)$.  Given the
generators we've found for $G_s$ and $1+{\cal N}$, the proposition
follows now immediately. \qed
\medskip

	Our next task is to relate ${\rm Aut}_g(S)$ to the group
$\widetilde{\rm Aut}{}^0(X)$ of Definition 4.1.  Every $\phi \in {\rm
Aut}_g(S)$ gives an automorphism $\phi_* : \cx^{\Delta(1)} \to
\cx^{\Delta(1)}$ which sends $t = (t_\rho) \in \cx^{\Delta(1)}$ to
$\phi_*(t) = (\phi(x_\rho)(t))$ (the notation $\phi(x_\rho)(t)$ means
the polynomial $\phi(x_\rho)$ evaluated at the point $t$).  It is easy
to check that $\phi_*$ which commutes with the action of $G$.  Our
next step is to show that $\phi_*$ preserves $\cx^{\Delta(1)} - Z
\subset \cx^{\Delta(1)}$:

\proclaim Proposition 4.6.  Assume that $X$ is a complete simplicial
toric variety, and let $\phi_*:\cx^{\Delta(1)} \to \cx^{\Delta(1)}$ is
the automorphism determined by $\phi \in {\rm Aut}_g(S)$.  Then
$\phi_*$ preserves $\cx^{\Delta(1)} - \Sigma$, commutes with $G$, and
hence lies in $\widetilde{\rm Aut}{}^0(X)$.

\prf Every element of $(\cx^*)^{\Delta(1)}$ clearly preserves
$\cx^{\Delta(1)} - Z$, so that by Proposition~4.5, it suffices to
prove the proposition when $\phi = y_m(\lambda)$, where $m$ is a root
in $R(N,\Delta)$.  Also, it suffices to show that
$y_m(\lambda)_*(\cx^{\Delta(1)}-Z) \subset \cx^{\Delta(1)}-Z$.  To
prove this, suppose $m \leftrightarrow (x_\rho,x^D)$.  We will write a
point of $\cx^{\Delta(1)}$ as $t = (t_\rho,{\bf t})$, where ${\bf t}$
is a vector indexed by $\Delta(1) - \{\rho\}$.  Then $y_m(\lambda)_*$
is given by
\d
y_m(\lambda)_*(t_\rho,{\bf t}) = (t_\rho + \lambda {\bf t}^D,{\bf t})\ ,
\d
where ${\bf t}^D$ is the monomial $x^D$ evaluated at ${\bf t}$ (this
makes sense since $x_\rho \notdiv x^D$).

	As in the proof of Theorem 2.1, write $\cx^{\Delta(1)}-Z
= \cup_\sigma U_\sigma$, where $U_\sigma = \{t \in \cx^{\Delta(1)} :
t^{\hat\sigma} \ne 0\}$.  Recall that $x^{\hat\sigma} = \prod_{\rho
\notin \sigma(1)} x_\rho$.  Now let $t = (t_\rho,{\bf t}) \in
U_\sigma$, and set $u = y_m(\lambda)_*(t) = (t_\rho + \lambda {\bf
t}^D,{\bf t})$.  We will show that $u \in U_\tau$ for some $\tau
\in \Delta$ by an argument similar to Proposition 3.14 of [Oda].  When
$\rho \in \sigma(1)$, we have $u^{\hat\sigma} = t^{\hat\sigma} \ne 0$,
so that $u \in U_\sigma$.  So suppose that $\rho \notin \sigma(1)$.
Since $\langle m,n_\rho\rangle = 1$ and $\langle m,n_{\rho'}\rangle
\le 0$ for $\rho' \ne \rho$, it follows that
$\check\sigma\cap\{m\}^\perp$ is a face of $\sigma$ such that $\rho$
is on one side of $\{m\}^\perp$ and all of the other one-dimensional
cones of $\Delta$ are on the other side.  But $\Delta$ is complete, so
$\rho$ and $\hat\sigma\cap\{m\}^\perp$ must be faces of some cone of
$\Delta$ which is on the same side of $\{m\}^\perp$ as $\rho$.  Since
$\Delta$ is simplicial, it follows that $\tau = \rho +
\hat\sigma\cap\{m\}^\perp$ must be a cone of $\Delta$.  There are now
two cases to consider:
\smallskip

\noindent Case 1: If $t^{\hat\tau} \ne 0$, then since $\rho \in
\tau(1)$, we have $u^{\hat\tau} = t^{\hat\tau} \ne 0$.  Thus $u \in
U_\tau$.
\smallskip

\noindent Case 2: If $t^{\hat\tau} = 0$, then $t_{\rho_0} = 0$ for some
$\rho_0 \notin \tau(1)$.  Since $\tau = \rho +
\hat\sigma\cap\{m\}^\perp$, this means $\langle -m,n_{\rho_0}\rangle >
0$.  Then $\prod_{\rho' \ne \rho} t_{\rho'}^{\langle
-m,n_{\rho'}\rangle} = 0$, which implies
\d
u^{\hat\sigma} = \Bigl(t_\rho + \lambda \prod_{\rho' \ne
\rho} t_{\rho'}^{\langle -m,n_{\rho'}\rangle}\Bigr) \prod_{\rho'
\notin \sigma(1)\cup\{\rho\}} t_{\rho'} = \bigl(t_\rho + 0\bigr)\prod_{\rho'
\notin \sigma(1)\cup\{\rho\}} t_{\rho'} = t^{\hat\sigma} \ne 0\ .
\d
It follows that $u \in U_\sigma$.
\smallskip

	Thus $y_m(\lambda)_*$ preserves $\cx^{\Delta(1)} -Z$, and
the proposition is proved. \qed
\medskip

	We can now begin the proof of Theorem 4.2:
\medskip

\noindent {\bf Proof of Theorem 4.2}.  We first prove part (iii) of
the theorem, which asserts that ${\rm Aut}_g(S)$ is naturally
isomorphic to $\widetilde{\rm Aut}{}^0(X)$.  From Proposition 4.6, the
map sending $\phi$ to $\phi_*^{-1}$ defines a group homomorphism
\d
{\rm Aut}_g(S) \longrightarrow \widetilde{\rm Aut}{}^0(X)\ .
\leqno{(12)}
\d
This map is easily seen to be injective. To see that it is surjective,
suppose that $\psi : \cx^{\Delta(1)} -Z \to \cx^{\Delta(1)} -Z$ is
$G$-equivariant.  Then for any $\rho$ we get $x_\rho\circ\psi
:\cx^{\Delta(1)} -Z \to \cx$, and since $\Sigma$ has codimension
$\ge 2$ by Lemma 1.4, purity of the branch locus implies that
$x_\rho\circ\psi$ extends to a polynomial $u_\rho$ defined on
$\cx^{\Delta(1)}$.  Since $\psi$ commutes with the action of $G$, it
follows that $u_\rho$ is homogeneous and $\deg(u_\rho) =
\deg(x_\rho)$.  Then $\phi(x_\rho) = u_\rho$ defines a graded
homomorphism $\phi : S \to S$ with the property that $\phi_* =
\psi$.  To see that $\phi$ is an isomorphism, apply the above argument
to $\psi^{-1}$.  This shows that (12) is an isomorphism, and part
(iii) of Theorem 4.2 is done.

	We now turn to part (ii) of the theorem, which asserts that we
have an exact sequence
\d
1 \longrightarrow G \longrightarrow \widetilde{\rm Aut}(X)
\longrightarrow {\rm Aut}(X) \longrightarrow 1\ .
\d
It is obvious that every element of $G$ gives the identity
automorphism of $X$.  Conversely, suppose that $\phi \in
\widetilde{\rm Aut}(X)$ induces the identity automorphism on $X$.  We
need to show that $\phi \in G$.  Our assumption implies that $\phi(Y)
= Y$ for every $G$-invariant $Y \subset \cx^{\Delta(1)}-Z$, so that in
particular, $\phi(U_\rho) = U_\rho$ for every $\rho \in \Delta(1)$.
Since $U_\rho = {\rm Spec}(S_\rho)$, we get an induced map $\phi^* :
S_\rho \to S_\rho$, and hence $\phi^*(S) \subset \cap_\rho S_\rho$.
It is easy to see that $\cap_\rho S_\rho = S$, and it follows that
$\phi^*(S) \subset S$.  Thus, for every $\rho$, we have
$\phi^*(x_\rho) = A_\rho \in S$.  Note, however, that $A_\rho$ need
not be homogeneous of degree $\deg(x_\rho)$.

	Since ${\bf V}(x_\rho) - Z \subset \cx^{\Delta(1)} -Z$ is also
$G$-invariant, it is stable under $\phi$.  But the $\rho$th coordinate
of $\phi(t)$ is $A_\rho(t)$, so that $A_\rho$ vanishes on ${\bf
V}(x_\rho) - Z$.  Since $Z$ has codimension $\ge 2$, we conclude that
$A_\rho$ vanishes on ${\bf V}(x_\rho)$.  The Nullstellensatz implies
that $x_\rho | A_\rho$, and hence $\phi^*(x_\rho) = x_\rho B_\rho$ for
some $B_\rho \in S$.  The same argument applied to $\phi^{-1}$ shows
that $\phi^{-1*}(x_\rho) = x_\rho C_\rho$ for some $C_\rho \in S$.
Then
\d
x_\rho = \phi^*\circ\phi^{-1*}(x_\rho) = \phi^*(x_\rho C_\rho) =
x_\rho B_\rho \phi^*(C_\rho)\ .
\d
This implies that $1 = B_\rho \phi^*(C_\rho)$ in the polynomial ring
$S=\cx[x_\rho]$, so that $B_\rho$ is a nonzero constant.  Thus $\phi =
(B_\rho) \in (\cx^*)^{\Delta(1)}$.  Since $\phi$ induces the identity
on $T \subset X$, the exact sequence $1 \to G \to (\cx^*)^{\Delta(1)}
\to T \to 1$ shows that $\phi \in G$ as claimed.  This proves that $G$
is the kernel of $\widetilde{\rm Aut}(X) \to {\rm Aut}(X)$.

	We still need to prove that $\widetilde{\rm Aut}(X) \to {\rm
Aut}(X)$ is onto.  We will first show that there is an exact
sequence
\d
1 \longrightarrow G \longrightarrow \widetilde{\rm Aut}{}^0(X)
\longrightarrow {\rm Aut}(X) \longrightarrow {\rm Aut}(A_{n-1}(X))
\ .\leqno{(13)}
\d
The map ${\rm Aut}(X) \to {\rm Aut}(A_{n-1}(X))$ is induced by direct
image of divisors: given $\phi \in {\rm Aut}(X)$, a divisor class
$[\sum_i a_i D_i]$ maps to $[\sum_i a_i\phi(D)]$.

	To show exactness, let $\psi \in \widetilde{\rm Aut}{}^0(X)$.
We need to show that the induced automorphism of $X$, which will also
be denoted $\psi$, is the identity on $A_{n-1}(X)$.  By part (iii) of
the theorem, we can write $\psi(t) = (u_\rho(t))$ where $u_\rho \in
S_{\deg(x_\rho)}$.  If $D'_\rho = {\bf V}(x_\rho) -Z$, then
\d
{\rm div}(u_\rho/x_\rho) = \psi^{-1}(D'_\rho) - D'_\rho
\d
in $\cx^{\Delta(1)}-Z$.  Since $u_\rho$ and $x_\rho$ have the
same degree, the quotient $u_\rho/x_\rho$ is $G$-invariant and
descends to a rational function on $X$.  Then it follows immediately
that
\d
{\rm div}(u_\rho/x_\rho) = \psi^{-1}(D_\rho) - D_\rho
\d
in $X$.  This shows that $[\psi^{-1}(D_\rho)] = [D_\rho]$, so that
$\psi$ induces the identity on $A_{n-1}(X)$.

	Next, assume that $\phi \in {\rm Aut}(X)$ induces
the identity on $A_{n-1}(X)$.  We need to show that $\phi$ comes from
an element of $\widetilde{\rm Aut}{}^0(X)$.  Fix $\rho \in \Delta(1)$,
and let $\a_i = \deg(x_\rho) = [D_\rho]$.  Since $[\phi^{-1}(D_\rho)] =
[D_\rho]$ in $A_{n-1}(X)$, we can find a rational function $g_\rho$ on
$X$ such that
\d
{\rm div}(g_\rho) = \phi^{-1}(D_\rho) - D_\rho\ .
\d
Then ${\rm div}(g_\rho) + D_\rho \ge 0$, which implies $g_\rho \in
H^0(X,{\cal O}_X(D_\rho))$.  Using the isomorphism $H^0(X,{\cal
O}_X(D_\rho)) \simeq S_{\a_i}$ from Proposition 1.1, we can write
$g_\rho = u_\rho/x_\rho$ for some $u_\rho \in S_{\a_i}$.  Note that
$u_\rho$ is only determined up to a constant.  We will choose the
constant as follows.

	Given $m \in M$, $\prod_\rho x^{\langle m,n_\rho\rangle}_\rho$
and $\prod_\rho u^{\langle m,n_\rho\rangle}_\rho$ are rational
functions on $X$.  Then
\d
\eqalign{{\rm div}\bigl(\textstyle{\prod_\rho u^{\langle
	m,n_\rho\rangle}_\rho}\bigr) &=
	{\rm div}\bigl(\textstyle{\prod_\rho (u_\rho/x_\rho)^{\langle
	m,n_\rho\rangle}}\bigr) + {\rm div}\bigl(\textstyle{\prod_\rho
	x^{\langle m,n_\rho\rangle}_\rho}\bigr) \cr
&= \textstyle{\sum_\rho} \langle m,n_\rho\rangle (\phi^{-1}(D_\rho)-
	D_\rho)  + \textstyle{\sum_\rho} \langle m,n_\rho\rangle D_\rho\cr
&= \phi^{-1}\bigl(\textstyle{\sum_\rho} \langle m,n_\rho\rangle
	D_\rho\bigr) = {\rm div}\Bigl(\bigl(\textstyle{\prod_\rho x^{\langle
	m,n_\rho\rangle}_\rho}\bigr)\circ\phi\Bigr)\ .\cr}
\d
This implies that there is $\gamma(m) \in \cx^*$ such that
\d
\gamma(m)\textstyle{\prod_\rho u^{\langle m,n_\rho\rangle}_\rho} =
	\bigl(\textstyle{\prod_\rho x^{\langle
m,n_\rho\rangle}_\rho}\bigr)\circ \phi\ .
\d
It is easy to see that $\gamma : M \to \cx^*$ is a homomorphism, which
implies that $\gamma \in T$.  The map $(\cx^*)^{\Delta(1)} \to T$ is
surjective, so that there is $t = (t_\rho) \in (\cx^*)^{\Delta(1)}$
such that $\gamma(m) = \prod_\rho t^{\langle m,n_\rho\rangle}_\rho$
for all $m \in M$.  Thus, if we replace $u_\rho$ by $t_\rho u_\rho$,
then for $m \in M$, we have
\d
\textstyle{\prod_\rho u^{\langle m,n_\rho\rangle}_\rho} =
	\bigl(\textstyle{\prod_\rho x^{\langle
m,n_\rho\rangle}_\rho}\bigr)\circ \phi\ .\leqno{(14)}
\d

	Now define $\psi : S \to S$ by $\psi(x_\rho) = u_\rho$, and
consider the diagram:
\d
\matrix{S_{\a_i} &&\mapname{\psi} && S_{\a_i} \cr
\ \ \downarrow \scriptstyle{x_\rho^{-1}} &&&&\ \ \downarrow
\scriptstyle{x_\rho^{-1}}\cr
H^0(X,{\cal O}_X(D_\rho)) & \mapname{\phi^*} & H^0(X,{\cal
O}_X(\phi^{-1}(D_\rho))) & \mapname{u_\rho/x_\rho} & H^0(X,{\cal
O}_X(D_\rho))\cr}
\d
Here, the vertical arrows are the isomorphisms given by Proposition
1.1.  Also, $\phi^*$ maps $g \in H^0(X,{\cal O}_X(D_\rho))$ to $g\circ
\phi \in H^0(X,{\cal O}_X(\phi^{-1}(D_\rho)))$, and the final map
on the bottom takes $g \in H^0(X,{\cal O}_X(\phi^{-1}(D_\rho)))$ to $g
\cdot (u_\rho/x_\rho) \in H^0(X,{\cal O}_X(D_\rho))$.  To see that
the diagram commutes, let $u \in S_{\a_i}$ be a monomial.  Using the
left and bottom maps, we get
\d
{u\over x_\rho}\circ\phi \cdot {u_\rho\over x_\rho} \in H^0(X,{\cal
O}_X(D_\rho))\ .
\d
However, since $u$ and $x_\rho$ have the same degree, we can write
$u/x_\rho = \prod_\rho x^{\langle m,n_\rho\rangle}_\rho$ for some $m
\in M$.  Then (14) easily implies that $(u/x_\rho)\circ\phi =
\psi(u)/u_\rho$, so that
\d
{u\over x_\rho}\circ\phi \cdot {u_\rho\over x_\rho} = {\psi(u) \over
u_\rho} \cdot {u_\rho\over x_\rho} = {\psi(u) \over x_\rho}\ .
\d
This shows that the above diagram commutes.  Since the maps on the
bottom are clearly isomorphisms, it follows that $\psi_{|S_{\a_i}}$ is
an isomorphism.  This shows that $\psi \in {\rm Aut}_g(S)$, and then
Proposition 4.6 gives $\psi_* \in \widetilde{\rm Aut}{}^0(X)$.

	From $\psi_*$ we get at automorphism of $X$ which we will also
denote by $\psi_*$.  Since $u_\rho = x_\rho\circ\psi_*$ on
$\cx^{\Delta(1)} - Z$, it follows that for any $m \in M$, we have
\d
\textstyle{\prod_\rho u^{\langle m,n_\rho\rangle}_\rho} =
	\bigl(\textstyle{\prod_\rho x^{\langle
m,n_\rho\rangle}_\rho}\bigr)\circ \psi_*\ .
\d
{}From (14), we see that $\phi$ and $\psi_*$ agree as birational
automorphisms of $T$, which implies $\psi = \psi_*$ as automorphisms
of $X$.  This completes the proof that the sequence (13) is exact.

	The next step in showing that $\widetilde{\rm Aut}(X) \to {\rm
Aut}(X)$ is onto is to study the role played by the automorphisms of
the fan $\Delta$.  Let ${\rm Aut}(N,\Delta)$ be the group of lattice
isomorphisms $\varphi$ of $N$ which preserve the fan $\Delta$.  Such a
$\varphi$ permutes $\Delta(1)$ and this permutation determines
$\varphi$ uniquely since $\Delta(1)$ spans $N_\r$.  It follows that
${\rm Aut}(N,\Delta)$ is a finite group.

	We first relate ${\rm Aut}(N,\Delta)$ to $\widetilde{\rm
Aut}(X)$.  Given $\varphi \in {\rm Aut}(N,\Delta)$, the induced
permutation of $\Delta(1)$ can be regarded as a permutation matrix
$P_\varphi \in GL(\cx^{\Delta(1)})$.  We claim that $\varphi$ induces
an automorphism $\widehat\varphi^*:G \to G$ such that for any $g \in
G$, we have
\d
g \circ P_\varphi = P_\varphi \circ \widehat\varphi^*(g)\ .\leqno{(15)}
\d
To prove this, first note that $\varphi$ induces a dual automorphism
$\varphi^*:M\to M$.  Consider the following diagram:
\d
\matrix{0 & \longrightarrow & M\hfill & \longrightarrow & \z^{\Delta(1)}
	& \longrightarrow & A_{n-1}(X) & \longrightarrow & 0\cr
& & \downarrow\,\scriptstyle{\varphi^{*-1}} & & \downarrow
	\,\scriptstyle{P_\varphi}& & & &\cr
0 & \longrightarrow & M\hfill & \longrightarrow & \z^{\Delta(1)} &
	\longrightarrow & A_{n-1}(X) & \longrightarrow &0\cr}
\d
It is straightforward to check that the diagram commutes, which shows
that we get an induced isomorphism $\widehat\varphi : A_{n-1}(X) \to
A_{n-1}(X)$.  This in turn gives $\widehat\varphi^* : G \to G$ since
$G = {\rm Hom}_\z(A_{n-1}(X),\cx^*)$, and (15) now follows easily.

	From (15), we see that $P_\varphi$ normalizes $G$, and notice
also that $P_\varphi$ preserves $Z$ because $\varphi$ preserves
$\Delta$.  Hence $P_\varphi \in \widetilde{\rm Aut}(X)$ for all
$\varphi \in {\rm Aut}(N,\Delta)$.  Furthermore, it is easy to check
that the induced automorphism of $X$ is the usual way that $\varphi
\in {\rm Aut}(N,\Delta)$ acts on $X$.  By abuse of notation, we will
write $P_\varphi \in {\rm Aut}(X)$.

	Let $H \subset {\rm Aut}(X)$ be the image of $\widetilde{\rm
Aut}{}^0(X) \to {\rm Aut}(X)$.  By (13) and Proposition 4.3, it
follows that $H \simeq \widetilde{\rm Aut}{}^0(X)/G \simeq {\rm
Aut}_g(S)/G$ is a connected affine algebraic group with
$(\cx^*)^{\Delta(1)}/G \simeq T$ as maximal torus.  Also, (13) implies
that $H$ is normal in ${\rm Aut}(X)$.

	Now we can finally show that $\widetilde{\rm Aut}(X) \to {\rm
Aut}(X)$ is onto.  Suppose that $\phi \in {\rm Aut}(X)$.  Then we have
$\phi T \phi^{-1} \subset \phi H \phi^{-1} = H$, so that $\phi T
\phi^{-1}$ is a maximal torus in $H$.  Since all maximal tori of $H$
are conjugate, there is some $\psi \in H$ such that $\phi T \phi^{-1}
= \psi T \psi^{-1}$.  If we set $\mu = \psi^{-1}\circ \phi$, it
follows that $\mu T \mu^{-1} = T$ in ${\rm Aut}(X)$.  This easily
implies that $\mu(T) = T$ in $X$.  Thus $\mu_{|T}$ is an automorphism
of $T$.  But the automorphism group of a torus is generated by
translations and group automorphisms.  Thus $\mu_{|T} = g \circ t$,
where $g$ is a group automorphism and $t \in T$.  We claim that
$\mu\circ t^{-1} : X \to X$ is $T$-equivariant.  This is certainly
true on $T$ since $(\mu\circ t^{-1})_{|T} = g$, and the claim follows
since $T$ is dense in $X$.  By Theorem 1.13 of [Oda], this implies
$\mu\circ t^{-1} = P_\varphi$ for some $\varphi \in {\rm
Aut}(N,\Delta)$.

	It follows that $\phi = \psi\circ\mu = \psi\circ
P_\varphi\circ t$ in ${\rm Aut}(X)$.  Thus ${\rm Aut}(X)$ is generated
by $H$ and the $P_\varphi$ for $\varphi \in {\rm Aut}(N,\Delta)$.
Since every $P_\varphi$ comes from an element of $\widetilde{\rm
Aut}(X)$, $\widetilde{\rm Aut}(X) \to {\rm Aut}(X)$ is surjective.
This completes the proof of part (ii) of Theorem 4.2.

	Finally, consider part (i) of the theorem, which asserts that
$\widetilde{\rm Aut}(X)$ is an affine algebraic group with
$\widetilde{\rm Aut}{}^0(X)$ as the connected component of the
identity.  To show this, note that since ${\rm Aut}(X)$ is generated
by the image of $\widetilde{\rm Aut}{}^0(X)$ and the $P_\varphi$, it
follows from part (ii) of the theorem that $\widetilde{\rm Aut}(X)$ is
generated by $\widetilde{\rm Aut}{}^0(X)$ and $P_\varphi$ for $\varphi
\in {\rm Aut}(N,\Delta)$.  Since ${\rm Aut}(N,\Delta)$ is finite,
$\widetilde{\rm Aut}{}^0(X)$ has finite index in $\widetilde{\rm
Aut}(X)$.  We already know that $\widetilde{\rm Aut}{}^0(X)$ is a
connected affine algebraic group, and hence $\widetilde{\rm Aut}(X)$
is an affine algebraic group with $\widetilde{\rm Aut}{}^0(X)$ as the
connected component of the identity.  Theorem 4.2 is now proved.\qed
\medskip

	From the above proof, we can extract some additional
information about $\widetilde{\rm Aut}(X)$ and ${\rm Aut}(X)$.  For
example, $\widetilde{\rm Aut}(X)$ is generated by
$(\cx^*)^{\Delta(1)}$, the $y_m(\lambda)$ for $m \in R(N,\Delta)$, and
the $P_\varphi$ for $\varphi \in {\rm Aut}(N,\Delta)$.  To discuss
${\rm Aut}(X)$ in more detail, we need some notation: for a root $m
\in R(N,\Delta)$, we have $y_m(\lambda)_* \in \widetilde{\rm
Aut}{}^0(X)$, and we let $x_m(\lambda)$ denote the
resulting automorphism of $X$.  It is easy to check that this is the
same automorphism $x_m(\lambda)$ as defined by [Demazure, \S2.3] or
[Oda, \S3.4].  Then we can describe ${\rm Aut}(X)$ as follows:

\proclaim Corollary 4.7. Let $X$ be a complete simplicial toric
variety.  Then:
\vskip0pt
\item{(i)} ${\rm Aut}(X)$ is an affine algebraic group with $T$ as
maximal torus.
\item{(ii)} ${\rm Aut}(X)$ is generated by $T$, the $x_m(\lambda)$ for
$m \in R(N,\Delta)$, and the $P_\varphi$ for $\varphi \in {\rm
Aut}(N,\Delta)$.
\item{(iii)} If ${\rm Aut}^0(X)$ is the connected component of the
identity of ${\rm Aut}(X)$, then
\d
{\rm Aut}^0(X) \simeq \widetilde{\rm Aut}{}^0(X)/G \simeq {\rm
Aut}_g(S)/G\ .
\d
\item{(iv)} ${\rm Aut}^0(X)$ is the semidirect product of its
unipotent radical with a reductive group $G_s$ whose root system
consists of components of type ${\bf A}$.
\item{(v)} From the partition $\Delta(1) = \Delta_1\cup\cdots\cup
\Delta_s$ (where the $x_\rho$ have the same degree for $\rho \in
\Delta_i$), let $\Sigma_{\Delta_i}$ denotes the symmetric group on
$\Delta_i$.  Then $\prod_{i=1}^s \Sigma_{\Delta_i}$ is the Weyl group
of $G_s$.  Furthermore, we can regard this group as a normal subgroup
of ${\rm Aut}(N,\Delta)$, and there are natural isomorphisms
\d
{\rm Aut}(X)/{\rm Aut}^0(X) \simeq \widetilde{\rm
Aut}(X)/\widetilde{\rm Aut}{}^0(X) \simeq {\rm Aut}(N,\Delta)/\bigl(
\textstyle{\prod_{i=1}^s} \Sigma_{\Delta_i}\bigr)\ .
\d
\vskip0pt

\prf Parts (i) -- (iv) of the corollary follow easily from the proof
of Theorem 4.2 and the other results of this section.  As for part
(v), it is clear that we have an isomorphism ${\rm Aut}(X)/{\rm
Aut}^0(X) \simeq \widetilde{\rm Aut}(X)/\widetilde{\rm Aut}{}^0(X)$.
Furthermore, we can regard ${\rm Aut}(N,\Delta)$ as a subgroup of
$\widetilde{\rm Aut}(X)$, and then we have an isomorphism
\d
\widetilde{\rm Aut}(X)/\widetilde{\rm Aut}{}^0(X) \simeq {\rm
Aut}(N,\Delta)/({\rm Aut}(N,\Delta)\cap\widetilde{\rm Aut}{}^0(X))\ .
\d
Thus, regarding elements of $\prod_{i=1}^s \Sigma_{\Delta_i}$ as
permutation matrices, it suffices to prove that
\d
{\rm Aut}(N,\Delta)\cap\widetilde{\rm Aut}{}^0(X) =
\textstyle{\prod_{i=1}^s} \Sigma_{\Delta_i}\ .
\d
One direction is easy, for if $P_\varphi \in {\rm
Aut}(N,\Delta)\cap\widetilde{\rm Aut}{}^0(X)$, then by (13), we know
that $P_\varphi$ induces the identity automorphism on $A_{n-1}(X)$.
This shows that as a permutation of the $\rho's$, $P_\varphi$
preserves the degree of $x_\rho$, and it follows immediately that
$P_\varphi \in \prod_{i=1}^s \Sigma_{\Delta_i}$.

	Going the other way, suppose that $P \in \prod_{i=1}^s
\Sigma_{\Delta_i}$.  Since $\prod_{i=1}^s \Sigma_{\Delta_i} \subset
\prod_{i=1}^s GL(S_{\a_i}) \subset \widetilde{\rm Aut}{}^0(X)$, $P$
induces an automorphism $\phi$ of $X$.  On $\cx^{\Delta(1)} -Z$,
$P$ permutes the variables, so that on $X$, $\phi$ permutes the
$D_\rho$.  Since $T = X - \cup_\rho D_\rho$, it follows that $\phi(T)
= T$.  As we saw in the proof of Theorem 4.2, this implies that $\phi
= P_\varphi\circ t$ for some $\varphi \in {\rm Aut}(N,\Delta)$ and $t
\in T$.  Back in $\widetilde{\rm Aut}(X)$, we then have $P =
P_\varphi\circ u$ for some $u \in (\cx^*)^{\Delta(1)}$.  Since $P$ and
$P_\varphi$ are permutation matrices, we conclude that $P= P_\varphi$,
and it follows that $P \in {\rm Aut}(N,\Delta)\cap
\widetilde{\rm Aut}{}^0(X)$.  This completes the proof of the
corollary. \qed
\medskip

	When $X$ is smooth, this corollary is due to Demazure (see
[Demazure] or [Oda, \S3.4]).
\bigskip

\noindent{\bf References}
\bigskip

\ref M.~Audin, {\sl The Topology of Torus Actions on Symplectic
Manifolds\/}, Progress in Math. {\bf 93}, Birkh\"auser, Basel Boston
Berlin, 1991.
\medskip

\ref M.~Demazure, {\it Sous-groupes alg\'ebriques de rang maximum du
groupe de Cremona\/}, Ann. Sci.~Ecole Norm.~Sup.~{\bf 3} (1970),
507--588.
\medskip

\ref V.~Danilov, {\it The geometry of toric varieties\/}, Russian
Math.~Surveys {\bf 33} (1978), 97--154.
\medskip

\ref I.~Dolgachev, {\it Weighted projective varieties\/}, in {\sl
Group Actions and Vector Fields\/}, edited by J.~Carrell, Lecture
Notes in Math.~{\bf 956}, Springer-Verlag, Berlin Heidelberg New York,
1982, 34--71.
\medskip

\ref J.~Fine, in preparation.
\medskip

\ref J.~Fogarty and D.~Mumford, {\sl Geometry Invariant Theory},
Second Edition, Springer-Verlag, Berlin Heidelberg New York, 1982.
\medskip

\ref W.~Fulton, {\sl Introduction to Toric Varieties\/}, Princeton
University Press, Princeton, 1993.
\medskip

\ref R.~Hartshorne, {\sl Algebraic Geometry\/}, Springer-Verlag,
Berlin Heidelberg New York, 1977.
\medskip

\ref I.~Musson, {\it Differential operators on toric varieties\/},
preprint.
\medskip

\ref T.~Oda, {\sl Convex Bodies and Algebraic Geometry\/},
Springer-Verlag, Berlin Heidelberg New York, 1988.
\medskip

\ref J.-P.~Serre, {\it Faisceaux alg\'ebriques coh\'erents\/}, Ann.~of
Math.~{\bf 61} (1955), 197--278.
\medskip

\end

\title{Erratum to ``The Homogeneous Coordinate Ring of a Toric Variety''}

\author{David A. Cox}

\address{Department of Mathematics, Amherst College, Amherst, MA 01002}

\email{dac@math.amherst.edu}

\maketitle

\setcounter{section}{4}
\setcounter{theorem}{2}

My paper ``The Homogeneous Coordinate Ring of a Toric Variety'' [J.\ Algebraic Geometry {\bf 4} (1995), 17--50] has some incorrect statements before and during the proof of Proposition 4.3.  The purpose of this note is to correct these errors and give a valid proof of the proposition.  I am very grateful to Alexander Duncan for bringing this matter to my attention.

We will use the same notation as the paper, hereinafter referred to as \cite{c}.  The major error in the paper occurs in the discussion following the statement of Theorem 4.2 of \cite{c}, where we assert that the set $\mathrm{End}_g(S)$ of graded $\C$-algebra homomorphisms $\phi : S \to S$ with $\phi(1) = 1$ is a $\C$-algebra.  Here, $S$ is the homogeneous coordinate ring of the toric variety $X$. Nowadays $S$ is called the \emph{total coordinate ring} (see \cite{cls}) or the \emph{Cox ring} (see \cite{hk}).   
However, as pointed out to me by Duncan, the composition of two $\C$-algebra homomorphisms is again a $\C$-algebra homomorphism, but the same is not true for their sum.  So right now, the best we can say is that $\mathrm{End}_g(S)$ is a monoid under composition.  

The proof of Proposition 4.3 in \cite{c} was based on the faulty assumption that $\Endg(S)$ is a $\C$-algebra.  Hence the main task of this note is to give a correct proof of the proposition.  Here is the proposition, with some improvements suggested by the referee.

\begin{proposition}
Let $X$ be a complete toric variety, and let $S$ be its homogeneous coordinate ring.  Then
\begin{itemize}
\item[{\rm (i)}] $\mathrm{Aut}_g(S)$ is a connected affine algebraic group of dimension equal to $\sum_{i=1}^s |\Delta_i|\mathrm{dim}_\C\, S_{\alpha_i}$, and $(\C^*)^{\Delta(1)} \subset \mathrm{Aut}_g(S)$ is a maximal torus.
\item[{\rm (ii)}] The unipotent radical $R_u$ of $\aut(S)$ is isomorphic as a variety to an affine space of dimension $\sum_{i=1}^s |\Delta_i|(\mathrm{dim}_\C\, S_{\alpha_i} - |\Delta_i|)$.
\item[{\rm (iii)}] $\aut(S)$ has a closed subgroup $G_s$ isomorphic to
  the reductive group $\prod_{i=1}^s \mathrm{GL}(S_{\alpha_i}')$ of
  dimension $\sum_{i=1}^s |\Delta_i|^2$.  Also,
  $(\C^*)^{\Delta(1)} \subseteq G_s$.  
\item[{\rm (iv)}] $\aut(S)$ is isomorphic to the semidirect product $R_u \rtimes G_s$.
\end{itemize}
\end{proposition}

\begin{proof}
To simplify notation, we write the direct sum decompostion $ S_{\a_i} = S_{\a_i}' \oplus S_{\a_i}''$ from (7) of \cite{c} as $S_i = S_i' \oplus S_i''$.
Since elements of $\Endg(S)$ are $\C$-linear, preserve degrees, and are determined uniquely by their values on the variables $x_\rho$, we have a bijection of sets
\[
\Endg(S) \simeq \prod_{i=1}^s\homc(S_{i}',S_{i}),\leqno{(\mye1)}
\]
and we also have an injection
\[
\Endg(S) \hookrightarrow \prod_{i=1}^s\Endc(S_{i}).\leqno{(\mye2)}
\]
that is compatible with composition. 

We first show that $\mathrm{Im}(\Endg(S)) \subset \prod_{i=1}^s\Endc(S_{i})$ is a variety.  Recall from \cite{c} that $\phi(S_{i}'') \subset S_{i}''$.  It follows $\phi \in \Endc(S)$ corresponds via (\mye2) to a collection of matrices
\[
\begin{pmatrix} A_i & 0 \\ B_i & C_i \end{pmatrix} \in \Endc(S_{i}), \quad i = 1,\dots, s,\leqno{(\mye3)}
\]
where we use the canonical basis of $S_{i} = S_i'\oplus S_i''$ given by monomials of degree $\a_i$ to identify matrices with linear maps.  Note that the $S_{i}'$-columns $\big(\begin{smallmatrix} A_i \\ B_i\end{smallmatrix}\big)$ of (\mye3) are the data that make up the map (\mye1).   The matrices $C_i$ come from evaluating $\phi$ at monomials in $S_i''$, which are products of $\ge 2$ variables that lie in various $S_j$ for $j \ne i$ (this follows from $S_0 = \C$).  Hence the entries in $C_i$ are detemined by the matrices $A_j,B_j$ for $j \ne i$.  We will say more about this below.

One fact we will need is how (\mye1) relates to composition.  Suppose that $\phi, \psi \in \Endg(S)$ map to matrices
\[
\begin{pmatrix} A_i & 0 \\ B_i & C_i \end{pmatrix}, \ 
\begin{pmatrix} A_i' & 0 \\ B_i' & C_i' \end{pmatrix},\quad i=1,...,s.
\]
Since (\mye2) is compatible with composition, we see that $\phi\circ\psi$ corresponds to the products
\[
\begin{pmatrix} A_i & 0 \\ B_i & C_i \end{pmatrix} 
\begin{pmatrix} A_i' & 0 \\ B_i' & C_i' \end{pmatrix}
=
\begin{pmatrix} A_i A_i'& 0 \\ B_iA_i' + C_i B_i' & C_iC_i' \end{pmatrix} ,
\quad i=1,...,s.
\]
It follows that in the bijection (\mye1), we have
\[
\text{if } \phi \longleftrightarrow \begin{pmatrix} A_i\\B_i\end{pmatrix} \text{ and } \psi \longleftrightarrow \begin{pmatrix} A_i'\\B_i'\end{pmatrix}, \text{ then }
\phi\circ\psi \longleftrightarrow \begin{pmatrix} A_iA_i'\\B_iA_i'+C_iB_i'\end{pmatrix}. \leqno{(\mye4)} 
\]
This will be useful later in the proof.

The next step is to write down the equations that define $\Endg(S)$ inside of $\prod_{i=1} \Endc(S_{i})$ in (\mye2).  Our treatment is inspired by \cite[Prop.\ 5.12]{bg}.  The equations come from two sources:
\begin{itemize}
\item First, all of the entries in the upper right-hand block must be zero.  This is the ``$0$'' in (\mye3).  
\item Second, suppose that we have monomials $x^D, x^E \in S_{i}''$.  Given $\phi \in \Endg(S)$, we have 
\[
\phi(x^D) = \cdots + c_{ED}^i x^E + \cdots,
\]
where $c_{ED}^i$ is the corresponding entry in $C_i$ in (\mye3) for $\phi$.  But $x^D$ is a product of variables $x_{\rho_1} \cdots x_{\rho_\ell}$,  where we allow duplications.  Note that $x_{\rho_j} \notin S_{\a_i}'$ since $S_0 = \C$.  It follows that 
\begin{align*}
c_{ED}^i &= \text{coefficient of $x^E$ in } \phi(x^D)\\
            &= \text{coefficient of $x^E$ in  } \phi(x_{\rho_1})\cdots \phi(x_{\rho_\ell}).
\end{align*}
Each $\phi(x_{\rho_j})$ is a linear combination of monomials whose
coefficients are the corresponding entries in the matrices $A_{k_j},B_{k_j}$, where $x_{\rho_j} \in S_{k_j}$, i.e., $\deg(x_{\rho_j}) = \alpha_{k_j}$.  Hence we get an equation linking $c_{ED}^i$ with entries in $A_{k_j},B_{k_j}$, $j=1,\dots,\ell$.  
\end{itemize}

This analysis shows that $\Endg(S)$ is a linear algebraic monoid in the sense of \cite{p}.  Since $\aut(S)$ is the group of invertible elements of $\Endg(S)$, it follows from \cite{p} that $\aut(S)$ is an algebraic group.  

We will need the following characterization of which elements of $\Endg(S)$ are invertible:\ if $\phi \in \Endg(S)$ corresponds to matrices (\mye3),  then
\[
\phi \in \aut(S) \iff A_i, C_i \text{ are invertible for } i= 1,\dots,s.\leqno{(\mye5)}
\]
One direction is obvious.  For the other, suppose that the $A_i,C_i$ are all invertible.  Then consider the element $\psi \in \Endg(S)$ such that
\[
\psi \longleftrightarrow \begin{pmatrix} A_i^{-1}\\- C_i^{-1}B_iA_i^{-1}\end{pmatrix}
\]
via (\mye1).  
Using (\mye4), one obtains $\phi \circ \psi \longleftrightarrow \big(\begin{smallmatrix} I \\ 0 \end{smallmatrix}\big)$, so that $\phi\circ\psi$ is the identity.  But then the matrices associated to $\phi$ and $\psi$ multiply to the identity in each $\Endc(S_i)$, which means that the same is true when we reverse the order.  Hence $\psi\circ\phi$ is also the identity, which proves that $\phi \in \aut(S)$.  This completes the proof of (\mye5).

As in \cite{c}, let
\[
\mathcal{N} = \prod_{i=1}^s \homc(S_{i}',S_{i}'').
\]
To define $1 + \mathcal{N} \subset \Endg(S)$, we have to be careful since endomorphisms cannot be added.  We let $1+\mathcal{N}$ consist of all
$\phi \longleftrightarrow \big(\begin{smallmatrix} I \\ B_i\end{smallmatrix}\big)$ via (\mye1), where $B_i \in \homc(S_i',S_i'')$.
Then an element $\phi \in 1 + \mathcal{N}$  gives matrices
\[
\begin{pmatrix} I & 0 \\ B_i & C_i \end{pmatrix} \in \Endc(S_{i}) \quad i = 1,\dots, s.\leqno{(\mye6)}
\]
We claim that these matrices are all unipotent.   

To study $C_i$, we order the monomials in $S_{i}''$ so that $x^D$ appears before $x^E$ whenever the total degree of $x^D$ (as a monomial in the polynomial ring $S$) is strictly smaller than the total degree of $x^E$.  Take $x^D \in S_{i}$ and write $x^D = x_{\rho_1}\cdots x_{\rho_l}$, so that $x^D$ has total degree $\ell$.   Applying $\phi$, we get 
\[
\phi(x_\rho) = \prod_{j=1}^\ell \phi(x_{\rho_j}) = \prod_{j=1}^\ell \big(x_{\rho_j} +{\textstyle \sum_{x^E \in S''_{\mathrm{deg(x_{\rho_j})}}}}b_{E,j} x^E\big)
\]
Since every monomial in $S_{\mathrm{deg(x_{\rho_j})}}''$ has total degree at least two, multiplying out the last product gives
\[
\phi(x^D) = x^D  + \text{terms of higher total degree}.
\]
Given how the monomials in $S_{\a_i}''$ are ordered, it follows that $C_i$ is lower triangular with $1$'s on the main diagonal.  Then the same is true for (\mye6), so that (\mye6) is unipotent as claimed.

Now that we know that $C_i$ is invertible, (\mye5) and (\mye6) imply that $\phi$ is invertible.  Hence we have proved that  $1+ \mathcal{N} \subset \aut(S)$.  Notice also that $1+\mathcal{N}$ is a closed subgroup of $\aut(S)$ by (\mye4) and (\mye6).  

Now we get to our main task, which is to establish the exact sequence of groups
\[
1 \longrightarrow 1+\mathcal{N} \stackrel{\alpha}{\longrightarrow} \aut(S) \stackrel{\beta}{\longrightarrow} \prod_{i=1}^s \mathrm{GL}(S_i') \longrightarrow 1.\leqno{(\mye7)}
\]
This is the exact sequence (9) of \cite{c}.  

The map $\alpha$ is the inclusion $1+\mathcal{N} \subset \aut(S)$ proved above.  The map $\beta$ is also easy to describe:\ if $\phi \in \aut(S)$ is specified by $\big(\begin{smallmatrix} A_i \\ B_i \end{smallmatrix}\big)$, then the $A_i$ are invertible by (\mye5) and hence give an element of $\prod_{i=1}^s \mathrm{GL}(S_i')$.  This is $\beta(\phi)$.   Note that $\beta$ is a group homomorphism by (\mye4).

The map $\alpha$ is clearly injective, and (\mye7) is exact at $\aut(S)$ by the definition of $1+\mathcal{N}$.  It remains to prove that $\beta$ is onto.  Suppose that we have invertible matrices $A_i \in \mathrm{GL}(S_i')$ for $i = 1,\dots,s$.  Then consider $\phi,\psi \in \Endg(S)$ such that
\[
\phi  \longleftrightarrow \begin{pmatrix} A_i\\0\end{pmatrix} \text{ and } \psi \longleftrightarrow \begin{pmatrix} A_i^{-1}\\0\end{pmatrix}\leqno{(\mye8)}
\]
via (\mye1).  Using (\mye4), one computes that
\[
\phi\circ\psi \longleftrightarrow \begin{pmatrix} A_iA_i^{-1}\\  0 \cdot A_i^{-1}+C_i\cdot 0\end{pmatrix} = \begin{pmatrix} I \\ 0\end{pmatrix}.
\]
This proves that $\phi\circ \psi$ is the identity, and switching $\phi$ and $\psi$ shows that $\psi\circ \phi$ is the identity as well.  Thus $\phi \in \aut(S)$.  Since $\beta$ maps $\phi$ to the $A_i$, surjectivity follows.  

Hence (\mye7) is exact, and we also know that $1+\mathcal{N}$ is unipotent.  Then part (ii) of the proposition follows because, as a variety, we have $1+\mathcal{N} \simeq \mathcal{N}$, which is an affine space of the required dimension.

For part (iii), note that the first half of (\mye8) gives a section 
\[
s^* : \prod_{i=1}^s \mathrm{GL}(S_i') \longrightarrow \aut(S)
\]
of the exact sequence (\mye7).  Note that $s^*$ is a group homomorphism by (\mye4).  The image is easily seen to be an algebraic subgroup containing $(\C^*)^{\Delta(1)}$.  This proves part (iii) of the proposition, and parts (iv) and (i) now follow without difficulty in view of (\mye7).  The proof is complete.
\end{proof}

Here are some final comments:
\begin{itemize}
\item Lemma 1.3 of \cite{c} is only used in the invalid proof of Proposition 4.3 in \cite{c}.  Hence this lemma can be ignored when reading the paper.  In the proof of Proposition 4.3 presented here, the ordering of Lemma~1.3 is replaced by the total degree ordering on the polynomial ring $S$.
\item The sentence following the first display in the proof of Proposition~4.5 of \cite{c} needs to be modified:\ ``is the maximal torus and hence lies in $G_s$" should be ``is a maximal torus contained in $G_s$.''
\item In \cite{d}, Demazure gives a functorial construction of the automorphism group of a toric variety $X$.  In \cite{c}, the approach is more concrete, based on the construction of $\aut(S)$ as a matrix group.  It would be useful to show that these two methods lead to the same algebraic group. 
\end{itemize}

\noindent {\bf Acknowledgements.} Besides Alexander Duncan, I am grateful to the referee for the reference \cite{p} and other helpful suggestions.